\documentclass{article}
\usepackage{latexsym}
\usepackage{amssymb}
\usepackage{amsmath}
\usepackage{amscd}
\usepackage[only,twlrm,ninrm,sevrm]{rawfonts}
\usepackage[latin1]{inputenc}
\numberwithin{equation}{section} \textwidth=5.45in
\newtheorem{Th}{Theorem}[section]
\newtheorem{Prop}{Proposition}[section]
\newtheorem{Co}{Corollary}[section]
\newtheorem{Le}{Lemma}[section]
\newtheorem{Remark}{Remark}[section]

\newcommand{\R}{\mathrm{I\!R\!}}

\newcommand{\N}{\mathrm{I\!N\!}}

\newcommand{\var}{\varepsilon}
\newcommand{\dis}{\displaystyle}
\newcommand{\rig}{\rightarrow}
\newcommand{\cha}{\widehat}
\newcommand{\righ}{\rightharpoonup}
\newcommand{\vart}{\vartheta}
\newcommand{\fim}{\blacksquare}

\begin{document}
\title{Multiplicity and concentration behavior of positive solutions for a Schr\"{o}dinger-Kirchhoff type problem via penalization method}
\author{
Giovany M. Figueiredo\footnote{Partially supported by CNPq/PQ
301242/2011-9 and 200237/2012-8} and Jo\~ao R. Santos J\'unior\footnote{Partially supported by CAPES - Brazil - 7155123/2012-9}\\
Faculdade de Matem\'atica\\
Universidade Federal do Par\'a\\
66.075-110 - Bel\'em - Par\'a - Brazil\\}

\date{}
\maketitle \vspace*{-25pt}

%%%%%%%%%%%%%%%%%%%%%%%%%%%%%%%%%%%%%%%%%%%%%%%%%%%%%%%%%%%%%%%%%%%%%%%%%%%%%%%%%%%%%%%%%%%%%%%%%%%%%%%%%%%%%%%%%%%%%%%%
%%%%%%%%%%%%%%%%%%%%%%%%%%%%%%%%%%%%%%%%%%%%%%%%%%%%%%%%%%%%%%%%%%%%%%%%%%%%%%%%%%%%%%%%%%%%%%%%%%%%%%%%%%%%%%%%%%%%%%%%
%%%%%%%%%%%%%%%%%%%%%%%%%%%%%%%%%%%%%%%%%%%%%%%%%%%%%%%%%%%%%%%%%%%%%%%%%%%%%%%%%%%%%%%%%%%%%%%%%%%%%%%%%%%%%%%%%%%%%%%%

\begin{abstract}
In this paper we are concerned with questions of multiplicity and concentration behavior of positive solutions of the elliptic problem
$$
 \left\{
\begin{array}{rcl}
\mathcal{L}_{\var}u=f(u) \ \ \mbox{in} \ \ \R^{3},\\
u>0 \ \ \mbox{in} \ \ \R^{3},\\
u \in H^{1}(\R^3),
\end{array}
\right.\leqno{(P_{\var})}
$$
where $\var$ is a small positive parameter,  $f:\mathbb{R}\rightarrow \mathbb{R}$ is a continuous function, $\mathcal{L}_{\var}$ is a nonlocal operator defined by
$$
\mathcal{L}_{\var}u=M\left(\dis\frac{1}{\var}\int_{\R^{3}}|\nabla
u|^{2}+\frac{1}{\var^{3}}\dis\int_{\R^{3}}V(x)u^{2}\right)\left[-\var^{2}\Delta
u + V(x)u \right],
$$
$M:\R_{+}\rig \R_{+}$ and $V:\R^{3}\rig \R$ are continuous functions
which verify some hypotheses.

\vspace{0.5cm}


\end{abstract}

%%%%%%%%%%%%%%%%%%%%%%%%%%%%%%%%%%%%%%%%%%%%%%%%%%%%%%%%%%%%%%%%%%%%%%%%%%%%%%%%%%%%%%%%%%%%%%%%%%%%%%%%%%%%%%%%%%%%%%%%
%%%%%%%%%%%%%%%%%%%%%%%%%%%%%%%%%%%%%%%%%%%%%%%%%%%%%%%%%%%%%%%%%%%%%%%%%%%%%%%%%%%%%%%%%%%%%%%%%%%%%%%%%%%%%%%%%%%%%%%%
%%%%%%%%%%%%%%%%%%%%%%%%%%%%%%%%%%%%%%%%%%%%%%%%%%%%%%%%%%%%%%%%%%%%%%%%%%%%%%%%%%%%%%%%%%%%%%%%%%%%%%%%%%%%%%%%%%%%%%%%

\section{Introduction} In this paper we shall focus our attention
on questions of multiplicity, concentration behavior and positivity of solutions for the following
problem
$$
 \left\{
\begin{array}{rcl}
\mathcal{L}_{\var}u=f(u) \ \ \mbox{in} \ \ \R^{3},\\
u>0 \ \ \mbox{in} \ \ \R^{3},\\
u \in H^{1}(\R^3),
\end{array}
\right.\leqno{(P_{\var})}
$$
where $\var$ is a small positive parameter,  $f:\mathbb{R}\rightarrow \mathbb{R}$ is a continuous function, $\mathcal{L}_{\var}$ is a nonlocal operator
defined by
$$
\mathcal{L}_{\var}u=M\left(\dis\frac{1}{\var}\int_{\R^{3}}|\nabla
u|^{2}+\frac{1}{\var^{3}}\dis\int_{\R^{3}}
V(x)u^{2}\right)\left[-\var^{2}\Delta u + V(x)u \right],
$$
and $M:\R_{+}\rig \R_{+}$ and $V:\R^{3}\rig \R$ are continuous
functions that satisfy some conditions which will be stated later
on.

Problem $(P_{\var})$ is a natural extension of two classes of
problems very important in applications, namely, Kirchhoff problems
and Schr\"{o}dinger problems.\\

$a)$ When $\var=1$ and $V=0$ problem $(P_{\var})$ we are dealing
with problem
$$
-M\biggl(\displaystyle\int_{\mathbb{R}^{3}}|\nabla u|^{2} dx
\biggl)\Delta u= f(u) \ \ \mbox{in} \ \ \mathbb{R}^{3},
$$
which represents the stationary case of Kirchhoff model for small
transverse vibrations of an elastic string by considering the
effects of the changes in the length of the string during the
vibrations.

Since that the length of string is variable during the vibrations,
then the tension in the string changes with time and depends of the
$L^{2}$ norm of the gradient of the displacement $u$. More
precisely, we have
$$
M(t)=\biggl(\displaystyle\frac{P_{0}}{h}+\displaystyle\frac{E}{2L}t\biggl),
\ \ t>0, \leqno{(1.2)}
$$
where $L$ is the length of the string, $h$ is the area of
cross-section, $E$ is the Young modulus of the material and $P_{0}$
is the initial tension.

Moreover, problem $(P_{\var})$  is called nonlocal because of the
presence of the term $M\left(\dis\int_{\mathbb{R}^{3}}|\nabla
u|^{2}\right)$ which implies that the equation in $(P_{\var})$ is no
longer a pointwise identity. This phenomenon causes some
mathematical difficulties which makes the study of such a class of
problem particularly interesting.

The version of problem $(P_{\var})$ in bounded domain began to call
attention of several researchers mainly after the work of Lions
$\cite{lions}$, where a functional analysis approach was proposed to
attack it.

We have to point out that nonlocal problems also appear in other
fields as, for example, biological systems where $u$ describes a
process which depends on the average of itself (for example,
population density). See, for example, \cite{alvescorreama} and its
references.

The reader may consult $\cite{alvescorrea}$,
$\cite{alvescorreagio}$, $\cite{alvescorreama}$, $\cite{Anelo}$,
$\cite{Anelo1}$, $\cite{Figueiredo and Santos Junior}$, $\cite{ma}$
and the references therein, for more informations on nonlocal
problems.\\

$b)$ On the other hand, when $M=1$ we have the problem
$$
 \left\{
\begin{array}{rcccl}
-\var^{2}\Delta u + V(x)u &=& f(u) \ \ \mbox{in} \ \ \R^{3},\\
u>0 \ \ \mbox{in} \ \ \R^{3},\\
u \in H^{1}(\R^3),
\end{array}
\right.\leqno{(1.3)}
$$
which arise in different models, for example, it is related with the
existence of standing waves of the nonlinear Schrodinger equation
$$
i\var\frac{\partial \Psi}{\partial
t}=-\var\Delta\Psi+(V(x)+E)\Psi-f(\Psi), \ \forall x\in
\R^{N},\leqno{(1.4)}
$$
where $f(t)=|t|^{p-2}u$ and $2<p<2^{\ast}=\dis\frac{2N}{N-2}$. A
standing wave of $(1.4)$ is a solution of the form
$\Psi(x,t)=\dis\exp(-iEt/\var)u(x)$. In this case, $u$ is a solution
of $(1.3)$. Existence and concentration of positive solutions for
the problem $(1.3)$ have been extensively studied in recent years,
see for example the papers $\cite{Ambrosetti and Badiale}$,
$\cite{Ambrosetti and Malchiodi}$, $\cite{Cingolani}$, $\cite{Del
Pino}$, $\cite{Floer}$, $\cite{Rabinowitz}$ and their references.\\

A considerable effort has been devoted during the last years in
studying  problems of the type $(P_{\var})$, as can be seen in
\cite{Alves and Figueiredo}, \cite{He}, \cite{Li}, \cite{Nie},
\cite{Wang 1},  \cite{Wu} and references therein. This is due to
their significance in applications as well as to their mathematical relevance.\\

Before stating our main result, we need the following hypotheses on
the function $M$:

\begin{description}

\item[($M_{1}$)] There is $m_{0}>0$ such that $M(t)\geq m_{0}, \ \forall t\geq 0$.

\item[($M_{2}$)] The function $t\mapsto M(t)$ is increasing.

\item[($M_{3}$)] There is $\vart\geq m_{0}>0$ such that

$$
\frac{M(t_{1})}{t_{1}}-\frac{M(t_{2})}{t_{2}}\leq \vart\left(
\frac{1}{t_{1}}-\frac{1}{t_{2}}\right),
$$
for all $t_{1}> t_{2}>0$.

\end{description}

The potential $V$ is a continuous function satisfying:
\begin{description}
\item[($V_{1}$)] There is $V_{0}>0$ such that $V_{0}=\dis\inf_{x\in\R^{3}}V(x)$.

\item[($V_{2}$)] For each $\delta>0$ there is a bounded and Lipschitz domain $\Omega \subset \R^{3}$ such that

$$
V_{0}<\dis \min_{\partial \Omega}V, \ \ \ \Pi=\{x\in \Omega: V(x)=V_{0}\}\neq \emptyset
$$
and
$$
\Pi_{\delta}=\{ x\in \R^{3}: dist(x, \Pi)\leq \delta\}\subset \Omega.
$$
\end{description}

Moreover, we assume that the continuous function $f$ vanishes in $(-\infty, 0)$ and verifies
\begin{description}
\item[($f_{1}$)] $$\lim_{t\to 0^{+}}\frac{f(t)}{t^{3}}=0. $$

\item[($f_{2}$)] There is $q\in (4,6)$ such that

$$
\lim_{t\to \infty}\frac{f(t)}{t^{q-1}}=0.
$$

\item[($f_{3}$)] There is $\theta\in(4,6)$ such that
$$
0<\theta F(t)\leq f(t)t, \ \forall t>0.
$$

\item[($f_{4}$)] The application $$t\mapsto\frac{f(t)}{t^{3}}$$ is non-decreasing in $(0,\infty)$.

\end{description}

The main result of this paper is:

\begin{Th}\label{Theorem 1.1}
Suppose that the function $M$ satisfies $(M_{1})-(M_{3})$, the
potential $V$ satisfies $(V_{1})-(V_{2})$ and the function $f$
satisfies $(f_{1})-(f_{4})$. Then, given $\delta>0$ there is
$\overline{\var}=\overline{\var}(\delta)>0$ such that the problem
$(P_{\var})$ has at least $Cat_{\Pi_{\delta}}(\Pi)$ positive
solutions, for all $\var\in (0,\overline{\var})$. Moreover, if
$u_{\var}$ denotes one of these positive solutions and
$\eta_{\var}\in R^{3}$ its global maximum, then
$$
\lim_{\var\to 0}V(\eta_{\var})=V_{0}. \label{Theorem 2}
$$
\end{Th}

A typical example of function verifying the assumptions
$(M_{1})-(M_{3})$ is given by $\dis M(t)=m_{0}+bt$, where $m_{0}>0$
and $b>0$. More generally, any function of the form $\dis
M(t)=m_{0}+bt+\dis\sum_{i=1}^{k}b_{i}t^{\gamma_{i}}$ with $b_{i}\geq
0$ and $\gamma_{i}\in (0,1)$ for all $i\in \{1,2,\ldots, k\}$
verifies the hypotheses $(M_{1})-(M_{3})$.

A typical example of function verifying the assumptions $(f_{1})-(f_{4})$ is given by
$\dis f(t)=\dis\sum_{i=1}^{k}c_{i}(t^{+})^{q_{i}-1}$ with $c_{i}\geq 0$ not all null and
$q_{i}\in [\theta, 6)$ for all $i\in \{1,2,\ldots, k\}$.

\vspace{.4cm}

\vspace{.4cm}

Recently some authors have considered problems of the type
$(P_{\var})$. For example, He and Zou $\cite{He}$, by using
Lusternik-Schnirelmann theory and minimax methods,  proved a result
of multiplicity and concentration behavior for the following equation
$$
 \left\{
\begin{array}{rcccl}
-(\var^{2}a+b\var\dis\int_{\R^{3}}|\nabla u|^{2})\Delta u + V(x)u &=& f(u)\ \ \mbox{in} \ \ \R^{3}\\
u > 0 \ \ \mbox{in}\ \  \R^{3},\\
u \in H^{1}(\R^3),
\end{array}
\right.\leqno{(1.5)}
$$
assuming, between others hypotheses, that $f\in C^{1}(\R)$ has a
subcritical growth 3-superlinear and the potential $V$ verifies a
assumption introduced by Rabinowitz $\cite{Rabinowitz}$, namely,
$$
V_{\infty}=\dis\liminf_{|x|\to
\infty}V(x)>V_{0}=\dis\inf_{\R^{3}}V(x)>0.\leqno{(R)}
$$

In $\cite{Wang 1}$, Wang, Tian, Xu and Zhang have considered the
problem
$$
 \left\{
\begin{array}{rcccl}
-(\var^{2}a+b\var\dis\int_{\R^{3}}|\nabla u|^{2})\Delta u + V(x)u &=& \lambda f(u)+ |u|^{4}u \ \ \mbox{in} \ \ \R^{3}\\
u > 0 \ \ \mbox{in} \ \ \R^{3},\\
u \in H^{1}(\R^3).\\
\end{array}
\right.\leqno{(1.6)}
$$
Assuming that $f$ is only continuous, has subcritical growth
3-superlinear and the potential verifies $(R)$, the authors
showed that $(1.6)$ has multiple positive solutions when $\lambda$
is sufficiently large, by using Lusternik-Schnirelmann theory,
minimax methods and a approach as in $\cite{Szulkin}$ (see also $\cite{Szulkin 1}$).

Other results for the problem Sch\"{o}dinger-Kirchhoff type can be
seen in \cite{Alves and Figueiredo}, \cite{Li}, \cite{Nie},
\cite{Wu} and references therein.

Motivated by results found in \cite{Alves and Figueiredo}, $\cite{Del Pino}$,
\cite{He} and \cite{Wang 1}, we study multiplicity via Lusternik-Schnirelmann
theory and concentration behavior of solutions for the problem
$(P_{\var})$. Here we use the hypotheses $(V_{1})-(V_{2})$ that were
first introduced by Del Pino and Felmer $\cite{Del Pino}$ for laplacian case.
For $p$-laplacian case, see $\cite{Alves and Figueiredo 1}$.

We emphasize that, at least in our knowledge, does not exist in the
literature actually available results involving problems
Schr\"{o}dinger-Kirchhoff  type, where the potential is like that
introduced by Del Pino and Felmer $\cite{Del Pino}$. This is a
difficulty that occurs, possibly by competition between the growth
of nonlocal term and the growth of nonlinearity.

Here, we use the same type of truncation explored in \cite{Del Pino},
however, we make a new approach and  some estimates are totally
different, for example, we show that solution of truncated problem
is solution of the original problem with distinct arguments.

Moreover,  we completed the results found in \cite{Alves and
Figueiredo}, \cite{He} and \cite{Wang 1} in the following sense:

$1$ - Since $M$ is a function more general than those in $\cite{He}$
and $\cite{Wang 1}$, we have a additional difficulty. In general,
the weak limit of the Palais-Smale sequences is not weak solution of
the autonomous problem. We overcome this difficulty with assumptions
different from those found in \cite{Alves and Figueiredo}.

$2$ - Since the function $f$ is only continuous, we cannot use
standard arguments on the Nehari manifold. Hence, our result is
similar then those found in \cite{Wang 1}. However, since the
hypotheses on function V are different, our arguments are completely
different. Moreover, our result is for all positive lambda.

\vspace{.3cm}

The paper is organized as follows. In the Section 2 we show that the
auxiliary problem has a positive solution and we introduce some tools
needed for the multiplicity result, namely, Lemma $\ref{lema3.2.2}$ and
Proposition $\ref{Proposition 4.1.2}$.  In the Section 3 we
study the autonomous problem associated. This study allows us
to show that the auxiliary problem has multiple solutions.   In the
section 4 we prove the main result using Moser iteration method
\cite{Moser}.

%%%%%%%%%%%%%%%%%%%%%%%%%%%%%%%%%%%%%%%%%%%%%%%%%%%%%%%%%%%%%%%%%%%%%%%%%%%%%%%%%%%%%%%%%%%%%%%%%%%%%%%%%%%%%%%%%%%%%%%%%%%%%%%%%%%%%%%%%%%%
%%%%%%%%%%%%%%%%%%%%%%%%%%%%%%%%%%%%%%%%%%%%%%%%%%%%%%%%%%%%%%%%%%%%%%%%%%%%%%%%%%%%%%%%%%%%%%%%%%%%%%%%%%%%%%%%%%%%%%%%%%%%%%%%%%%%%%%%%%%%
%%%%%%%%%%%%%%%%%%%%%%%%%%%%%%%%%%%%%%%%%%%%%%%%%%%%%%%%%%%%%%%%%%%%%%%%%%%%%%%%%%%%%%%%%%%%%%%%%%%%%%%%%%%%%%%%%%%%%%%%%%%%%%%%%%%%%%%%%%%%

\section{The auxiliary problem}

Considering the change of variable $x=\var z$ in $(P_{\var})$ we obtain
the modified problem

$$
 \left\{
\begin{array}{rcl}
\mathcal{\widetilde{L}_{\var}}u=f(u) \ \ \mbox{in} \ \ \R^{3},\\
u>0 \ \ \mbox{in} \ \ \R^{3},\\
u \in H^{1}(\R^3),
\end{array}
\right.\leqno{(\widetilde{P}_{\var})}
$$
where
$$
\mathcal{\widetilde{L}_{\var}}u=M\left(\dis\int_{\R^{3}}|\nabla
u|^{2}+\dis\int_{\R^{3}}V(\var x)u^{2}\right)\left[-\Delta u +
V(\var x)u \right],
$$
which is clearly equivalent to $(P_{\var})$.

Since $(f_{1})$ and $(f_{4})$ imply that
$$
\displaystyle\lim_{t\rightarrow 0^{+}}\frac{f(t)}{t}=0
$$
and since that
$$
t\mapsto\frac{f(t)}{t}
$$
is a application in $(0,\infty)$ which is increasing and unbounded,
we can adapt to our case the penalization method introduced by Del
Pino and Felmer $\cite{Del Pino}$.

Let $K>\dis\frac{2}{m_{0}}$, where $m_{0}$ is given in $(M_{1})$ and
$a>0$ such that $f(a)=\dis\frac{V_{0}}{K}a$. We define
$$\widetilde{f}(t)=
\left\{
\begin{array}{ccc}
   f(t) & \mbox{if} & t\leq a , \\
  \frac{V_{0}}{K}t & \mbox{if} & t>a
\end{array}
\right.
$$
and
$$
g(x,t)=\chi_{\Omega}(x)f(t)+(1-\chi_{\Omega}(x))\widetilde{f}(t),
$$
where $\chi$ is characteristic function of set $\Omega$. From
hypotheses $(f_{1})-(f_{4})$ we get that $g$ is a Carath\'{e}odory
function and the following conditions are observed:

\begin{description}
\item[($g_{1}$)] $$\lim_{t\to 0^{+}}\frac{g(x,t)}{t^{3}}=0,$$ uniformly in $x\in \R^{3}$.

\item[($g_{2}$)]

 $$\lim_{t\to \infty}\frac{g(x,t)}{t^{q-1}}=0, \ \mbox{uniformly in} \ x\in \mathbb{R}^{3},$$

\item[($g_{3}$)]

$(i)$ $$ 0\leq \theta G(x,t)<g(x,t)t, \ \forall x\in \Omega \
\mbox{and} \ \forall t>0$$ and

$(ii)$ $$ 0\leq 2 G(x,t)\leq g(x,t)t\leq \frac{1}{K}V(x)t^{2}, \ \forall x\in \R^{3}\backslash\Omega \ \mbox{and} \ \forall t>0.$$

\item[($g_{4}$)]  For each $x\in \Omega$, the application $t\mapsto \frac{g(x,t)}{t^{3}}$ is increasing in $(0,\infty)$ and for
each $x\in \R^{3}\backslash\Omega$, the application $t\mapsto
\frac{g(x,t)}{t^{3}}$ is increasing in $(0,a)$.

\end{description}

\noindent Moreover, from definition of $g$, we have $g(x,t)\leq
f(t)$, for all $t\in (0, +\infty)$ and for all $x\in \R^{3}$,
$g(x,t)=0$ for all $t\in (-\infty,0)$ and for all $x\in \R^{3}$.

Now we study the auxiliary problem

$$
 \left\{
\begin{array}{rcl}
\mathcal{\widetilde{L}_{\var}}u=g(\var x, u), \ \R^{3}\\
u>0, \ \R^{3}\\
u \in H^{1}(\R^3).
\end{array}
\right.\leqno{(P_{\var, A})}
$$
Observe that positive solutions of $(P_{\var,A})$ with $u(x)\leq a$
for each $x\in \R^{3}\backslash\Omega$ are also positive solutions
of $(\widetilde{P}_{\var})$.

We obtain solutions of $(P_{\var,A})$ as critical points of the
energy functional
$$
J_{\var}(u)=\frac{1}{2}\cha{M}\left(\dis\int_{\R^{3}}|\nabla u|^{2}+\dis\int_{\R^{3}}V(\var x) u^{2}\right)-\dis\int_{\R^{3}}G(\var x, u),
$$
where $\cha{M}(t)=\dis\int_{0}^{t}M(s)ds$ and $G(x,t)=\dis\int_{0}^{t}g(\var x,s)ds$, which is well defined on
the Hilbert space $H_{\var}$, given by
$$
H_{\var}=\{u\in H^{1}(\R^{3}): \dis\int_{\R^{3}}V(\var x) u^{2}<\infty\},
$$
provided of the inner product
$$
(u,v)_{\var}=\dis\int_{\R^{3}}\nabla u\nabla
v+\dis\int_{\R^{3}}V(\var x) u v.
$$
The norm induced by inner product is denoted by
$$
\| u\|_{\var}^{2}=\dis\int_{\R^{3}}|\nabla u|^{2}+\dis\int_{\R^{3}}V(\var x) u^{2}.
$$
Since $M$ and $f$ are  continuous we have that $J_{\var}\in
C^{1}(H_{\var},\R)$ and
$$
J'_{\var}(u)v=M(\|
u\|_{\var}^{2})(u,v)_{\var}-\dis\int_{\R^{3}}g(\var x,u)v, \ \forall
u,v\in H_{\var}.
$$
Now, we will fix some notations. We denote the Nehari manifold associated to $J_{\var}$ by
$$
\mathcal{N}_{\var}=\{u\in H_{\var}\backslash\{0\}: J'_{\var}(u)u=0\}.
$$
We denote by $H_{\var}^{+}$ the open subset of $H_{\var}$
given by
$$
H_{\var}^{+}=\{u\in H_{\var}: |supp (u^{+})\cap \Omega_{\var}|>0 \},
$$
and $S_{\var}^{+}=S_{\var}\cap H_{\var}^{+}$, where
$S_{\var}$ is the unit sphere of $H_{\var}$.

Note that $S_{\var}^{+}$ is a incomplete
$C^{1,1}$-manifold of codimension $1$, modeled on $H_{\var}$ and contained  in the open
$H_{\var}^{+}$. Thus, $H_{\var}=T_{u}S_{\var}^{+}\oplus
\R \ u$ for each $u\in S_{\var}^{+}$, where
$T_{u}S_{\var}^{+}=\{v\in H_{\var}:(u, v)_{\var}=0\}$.

We also define the set $\Omega_{\var}$ by
$$
\Omega_{\var}=\{x\in \R^{3}: \var x\in\Omega\}.
$$

Finally, we mean by weak solution of $(P_{\var, A})$ a function $u\in H_{\var}$ such that
$$
M(\| u\|_{\var}^{2})(u,v)_{\var}=\dis\int_{\R^{3}}g(\var x,u)v, \ \forall v\in H_{\var}.
$$
Therefore, critical points of $J_{\var}$ are weak solutions of $(P_{\var,A})$.

\begin{Le}\label{Lemma 4.1.1}
The functional $J_{\var}$ satisfies the following conditions:

a) There are $\alpha, \rho>0$ such that
$$
J_{\var}(u)\geq \alpha, \ \mbox{with} \ \|u\|_{\var}=\rho.
$$

b) There is $e\in H_{\var}\backslash B_{\rho}(0)$ with $J_{\var}(e)<0$.

\end{Le}

\noindent {\bf Proof.} The item a) follows directly from the hypotheses $(M_{1})$, $(g_{1})$ and $(g_{2})$.

On the other hand, it follows from $(M_{3})$ that there is
$\gamma_{1}>0$ such that $M(t)\leq \gamma_{1}(1+t)$ for all $t\geq
0$. So, for each $u\in H_{\var}^{+}$ and $t>0$ we have
\begin{eqnarray*}
J_{\var}(tu)&=&\frac{1}{2}\cha{M}(\|tu\|_{\var}^{2})-\int_{\R^{3}} G(\var x, tu)\\
&\leq&\frac{\gamma_{1}}{2}t^{2}\|u\|_{\var}^{2}+\frac{\gamma_{1}}{4}t^{4}\|u\|_{\var}^{4}-\int_{\Omega_{\var}}
G(\var x, tu).
\end{eqnarray*}
From $(g_{3})(i)$, we obtain $C_{1}, C_{2}>0$ such that
$$
J_{\var}(tu)\leq
\frac{\gamma_{1}}{2}t^{2}\|u\|_{\var}^{2}+\frac{\gamma_{1}}{4}t^{4}\|u\|_{\var}^{4}
-C_{1}t^{\theta}\int_{\Omega_{\var}}(u^{+})^{\theta}+C_{2}|supp(u^{+})\cap\Omega_{\var}|.
$$
Since $\theta\in (4,6)$ we conclude b).
$\fim$

Once $f$ and $M$ are only continuous the next two results are very
important, because allow us to overcome the non-differentiability
of $\mathcal{N}_{\var}$ (see Lemma \ref{lema3.2.2} $(A_{3})$ and Proposition \ref{Proposition 4.1.2}) and
the incompleteness of $S_{\var}^{+}$ (see  Lemma \ref{lema3.2.2} $(A_{4})$).

\begin{Le}\label{lema3.2.2}
Suppose that the function $M$ satisfies $(M_{1})-(M_{3})$,
the potential $V$ satisfies $(V_{1})-(V_{2})$ and the function
$f$ satisfies $(f_{1})-(f_{4})$. So:

\begin{description}
\item[($A_{1}$)] For each $u\in H_{\var}^{+}$, let
$h:\R_{+}\rig\R$ be defined by $h_{u}(t)=J_{\var}(tu)$. Then,
there is a unique $t_{u}>0$ such that $h_{u}'(t)>0$ in $(0,t_{u})$
and $h_{u}'(t)<0$ in $(t_{u}, \infty)$.

\item[($A_{2}$)] there is $\tau>0$ independent on $u$ such that $t_{u}\geq \tau$ for all $u\in S^{+}_{\var}$. Moreover,
for each compact set $\mathcal{W}\subset S^{+}_{\var}$ there is
$C_{\mathcal{W}}>0$ such that $t_{u}\leq C_{\mathcal{W}}$, for all
$u\in \mathcal{W}$.\label{Lemma 4.1.2}

\item[($A_{3}$)] The map
$\widehat{m}_{\var}:H_{\var}^{+}\rightarrow
\mathcal{N}_{\var}$ given by $\widehat{m}_{\var}(u)=t_{u}u$ is
continuous and
$m_{\var}:=\widehat{m}_{\var_{\bigl|S^{+}_{\var}}}$ is a
homeomorphism between $S^{+}_{\var}$ and $\mathcal{N}_{\var}$.
Moreover, $m_{\var}^{-1}(u)= \frac{u}{\|u\|_{\var}}$.

\item[($A_{4}$)] If there is a sequence $(u_{n})\subset
S^{+}_{\var}$ such that \mbox{dist}$(u_{n},\partial
H^{+}_{\var})\rightarrow 0$, then
$\|m_{\var}(u_{n})\|_{\var}\rightarrow \infty$ and
$I_{\var}(m_{\var}(u_{n}))\rightarrow \infty$.
\end{description}
\end{Le}

\noindent {\bf Proof.} For proving $(A_{1})$, it is sufficient to note that,
from the Lemma $\ref{Lemma 4.1.1}$,
$h_{u}(0)=0$, $h_{u}(t)>0$ when $t>0$ is small and $h_{u}(t)<0$
when $t>0$ is large. Since $h_{u}\in C^{1}(\R_{+},\R)$, there is
$t_{u}>0$ global maximum point of $h_{u}$ and $h_{u}'(t_{u})=0$.
Thus, $J'_{\var}(t_{u}u)(t_{u}u)=0$ and $t_{u}u\in
\mathcal{N}_{\var}$. We see that $t_{u}>0$ is the unique positive
number such that $h_{u}'(t_{u})=0$. Indeed, suppose by
contradiction that there are $t_{1}>t_{2}>0$ with
$h_{u}'(t_{1})=h_{u}'(t_{2})=0$. Then, for $i=1,2$
$$
t_{i}M(\|t_{i}u\|_{\var}^{2})\|u\|_{\var}^{2}=\int_{\R^{3}} g(\var
x, t_{i}u)u.
$$
So,
$$
\frac{M(\|t_{i}u\|_{\var}^{2})}{\|t_{i}u\|_{\var}^{2}}=\frac{1}{\|u\|_{\var}^{4}}\int_{\R^{3}}\left[\frac{g(\var x, t_{i}u)}{(t_{i}u)^{3}}\right]u^{4}.
$$
Therefore,
$$
\frac{M(\|t_{1}u\|_{\var}^{2})}{\|t_{1}u\|_{\var}^{2}}-\frac{M(\|t_{2}u\|_{\var}^{2})}{\|t_{2}u\|_{\var}^{2}}=
\frac{1}{\|u\|_{\var}^{4}}\int_{\R^{3}}\left[\frac{g(\var x,
t_{1}u)}{(t_{1}u)^{3}}-\frac{g(\var x,
t_{2}u)}{(t_{2}u)^{3}}\right]u^{4}.
$$
It follows from $(M_{3})$ and $(g_{4})$ that
\begin{eqnarray*}
\frac{\vart}{\|u\|_{\var}^{2}}\left( \frac{1}{t_{1}^{2}}-
\frac{1}{t_{2}^{2}}\right)&\geq&\frac{1}{\|u\|_{\var}^{4}}\int_{(\R^{3}\backslash\Omega_{\var})\cap\{t_{2}
u\leq a<t_{1}u\}}\left[\frac{g(\var x, t_{1}u)}{(t_{1}u)^{3}}-\frac{g(\var x, t_{2}u)}{(t_{2}u)^{3}}\right]u^{4}\\
&+&\frac{1}{\|u\|_{\var}^{4}}\int_{(\R^{3}\backslash\Omega_{\var})\cap\{a<t_{2}u\}}
\left[\frac{g(\var x, t_{1}u)}{(t_{1}u)^{3}}-\frac{g(\var x,
t_{2}u)}{(t_{2}u)^{3}}\right]u^{4}.
\end{eqnarray*}
By using the definition of $g$ we obtain
\begin{eqnarray*}
\frac{\vart}{\|u\|_{\var}^{2}}\left( \frac{1}{t_{1}^{2}}-\frac{1}{t_{2}^{2}}\right)&\geq&\frac{1}{\|u\|_{\var}^{4}}\int_{(\R^{3}\backslash\Omega_{\var})\cap\{t_{2}u\leq a<t_{1}u\}}\left[\frac{V_{0}}{K}\frac{1}{(t_{1}u)^{2}}
-\frac{f(t_{2}u)}{(t_{2}u)^{3}}\right]u^{4}\\
&+&\frac{1}{\|u\|_{\var}^{4}}\frac{1}{K}\left(\frac{1}{t_{1}^{2}}-\frac{1}{t_{2}^{2}}\right)\int_{(\R^{3}\backslash
\Omega_{\var})\cap\{a<t_{2}u\}}V_{0}u^{2}.
\end{eqnarray*}
Multiplying both sides by
$\frac{\|u\|_{\var}^{4}}{\left(\frac{1}{t_{1}^{2}}-\frac{1}{t_{2}^{2}}\right)}$
and using the hypothesis $t_{1}>t_{2}$, it follows that
\begin{eqnarray*}
\vart \|u\|_{\var}^{2}&\leq&
\frac{t_{1}^{2}t_{2}^{2}}{t_{2}^{2}-t_{1}^{2}}\int_{(\R^{3}\backslash\Omega_{\var})\cap\{t_{2}u\leq a<t_{1}u\}}\left[\frac{V_{0}}{K}\frac{1}{(t_{1}u)^{2}}
-\frac{f(t_{2}u)}{(t_{2}u)^{3}}\right]u^{4}\\
&+&\frac{1}{K}\int_{(\R^{3}\backslash\Omega_{\var})\cap\{a<t_{2}u\}}V_{0}u^{2}.
\end{eqnarray*}
Thereby,
\begin{eqnarray*}
\vart \|u\|_{\var}^{2}&\leq&-\left(\frac{t_{2}^{2}}{t_{1}^{2}-t_{2}^{2}}\right)\frac{1}{K} \int_{(\R^{3}\backslash\Omega_{\var})\cap\{t_{2}u\leq a<t_{1}u\}}V_{0}u^{2}\\
&+&
\left(\frac{t_{1}^{2}}{t_{1}^{2}-t_{2}^{2}}\right)\int_{(\R^{3}\backslash\Omega_{\var})\cap\{t_{2}u\leq
a<t_{1}u\}}\frac{f(t_{2}u)}{t_{2}u}u^{2}
+\frac{1}{K}\int_{(\R^{3}\backslash\Omega_{\var})\cap\{a<t_{2}u\}}V_{0}u^{2}.
\end{eqnarray*}
So,
$$
\vart \|u\|_{\var}^{2}\leq\frac{1}{K}\int_{\R^{3}\backslash\Omega_{\var}}V_{0}u^{2}\leq \frac{1}{K}\|u\|_{\var}^{2}.
$$

Since $u\neq 0$, we have that $\vart\leq \frac{1}{K}<m_{0}$, but
this is a contradiction. Thus, $(A_{1})$ is proved.

$(A_{2})$ Now, let $u\in S_{\var}^{+}$. From $(M_{1})$, $(g_{1})$, $(g_{2})$ and
from the Sobolev embeddings
$$
m_{0}t_{u}\leq M(t_{u}^{2})t_{u}=\dis\int_{\R^{3}} g(\var x,
t_{u}u)u\leq
\frac{\xi}{4}C_{1}t_{u}^{3}+\frac{C_{\xi}}{q}C_{2}t_{u}^{q-1},
$$
since $\xi>0$ is arbitrary, we obtain $\tau>0$ such that $t_{u}\geq \tau$.
Finally, if
$\mathcal{W}\subset S_{\var}^{+}$ is compact, suppose by contradiction
that there is $\{u_{n}\}\subset \mathcal{W}$ such that
$t_{n}=t_{u_{n}}\rig \infty$. Since $\mathcal{W}$ is compact, there
is $u\in \mathcal{W}$ with $u_{n}\rig u$ in $H_{\var}$. It follows
from the arguments used in the proof of item b) of the Lemma
$\ref{Lemma 4.1.1}$ that $J_{\var}(t_{n}u_{n})\rig -\infty$. On the
other hand, note that if $v\in \mathcal{N}_{\var}$, then by
$(g_{3})(i)$
\begin{eqnarray*}
J_{\var}(v)&=&J_{\var}(v)-\frac{1}{\theta}J_{\var}'(v)v\\
&\geq&\frac{1}{2}\cha{M}(\|v\|_{\var}^{2})-\frac{1}{\theta}M(\|v\|_{\var}^{2})\|v\|_{\var}^{2}
+\frac{1}{\theta}\int_{\R^{3}\setminus\Omega_{\var}}\left[g(\var x,v)v-\theta G(\var x,v)\right].
\end{eqnarray*}
From $(g_{3})(ii)$ we have
$$
J_{\var}(v)\geq \frac{1}{2}\cha{M}(\|v\|_{\var}^{2})-\frac{1}{\theta}M(\|v\|_{\var}^{2})\|v\|_{\var}^{2}-\left(\frac{\theta - 2}{2\theta}\right)\frac{1}{K}\int_{\R^{3}\setminus\Omega_{\var}}V(\var x)v^{2},
$$
and so
$$
J_{\var}(v)\geq \frac{1}{2}\cha{M}(\|v\|_{\var}^{2})-\frac{1}{\theta}M(\|v\|_{\var}^{2})\|v\|_{\var}^{2}-\left(\frac{\theta - 2}{2\theta}\right)\frac{1}{K}\|v\|_{\var}^{2}.
$$
By using the hypothesis $(M_{3})$, we derive $\cha{M}(t)\geq \dis\frac{[M(t)+\vart]}{2}t$, for all $t\geq 0$. Thence,
$$
J_{\var}(v)\geq \left(\frac{\theta-4}{4\theta}\right)M(\|v\|_{\var}^{2})\|v\|_{\var}^{2}+\frac{\vart}{4}\|v\|_{\var}^{2}-\left(\frac{\theta - 2}{2\theta}\right)\frac{1}{K}\|v\|_{\var}^{2}.
$$
From $\vart\geq m_{0}$ and $(M_{1})$, we conclude
$$
J_{\var}(v)\geq \left(\frac{\theta-2}{2\theta}\right)\left(m_{0}-\frac{1}{K}\right)\|v\|_{\var}^{2}.
$$
Once $\{t_{n}u_{n}\}\subset \mathcal{N}_{\var}$, we obtain
$$
\frac{1}{t_{n}^{2}}J_{\var}(t_{n}u_{n})\geq \left(\frac{\theta-2}{2\theta}\right)\left(m_{0}-\frac{1}{K}\right), \ \forall n\in \N.
$$
However, choosing sufficiently large values of $n$
$$
0\geq \left(\frac{\theta - 2}{2\theta}\right)\left( m_{0}-\frac{1}{K}\right)>0,
$$
a contradiction. Therefore $(A_{2})$ is true.

$(A_{3})$ Firstly we observe that $\cha{m}_{\var}$, $m_{\var}$ and $m_{\var}^{-1}$ are well defined. In fact, by $(A_{1})$,
for each $u\in H_{\var}^{+}$, there is a unique $m_{\var}(u)\in \mathcal{N}_{\var}$.
On the other hand, if $u\in \mathcal{N}_{\var}$ then $u\in H_{\var}^{+}$.
Otherwise, we have $|supp(u^{+})\cap \Omega_{\var}|=0$ and by $(g_{3})(ii)$
$$
0<M(\|u\|_{\var}^{2})\|u\|_{\var}^{2}=\int_{\R^{3}}g(\var x, u)u=
\int_{\R^{3}\backslash\Omega_{\var}}g(\var x, u^{+})u^{+}\leq
\frac{1}{K}\int_{\R^{3}\backslash\Omega_{\var}}V(\var x)u^{2}.
$$
Hence, from $(M_{1})$
$$
0<\left(m_{0}-\frac{1}{K}\right)\|u\|_{\var}^{2}\leq 0,
$$
a contradiction. Consequently $m^{-1}_{\var}(u)=\frac{u}{\|u\|_{\var}}\in S_{\var}^{+}$,
$m_{\var}^{-1}$ is well defined and it is a continuous function. Since,
$$
m_{\var}^{-1}(m_{\var}(u))=m_{\var}^{-1}(t_{u}u)=\frac{t_{u}u}{t_{u}\|u\|_{\var}}=u, \ \forall \ u\in S_{\var}^{+},
$$
we conclude that $m_{\var}$ is a bijection.  To show that $\cha{m}_{\var}:H_{\var}^{+}\rig \mathcal{N}_{\var}$
is continuous, let $\{u_{n}\}\subset H_{\var}^{+}$ and $u\in H_{\var}^{+}$ be
such that $u_{n}\rig u$ in $H_{\var}$. From $(A_{2})$, there is $t_{0}>0$ such that
$t_{u_{n}}\rig t_{0}$. Since, $t_{u_{n}}u_{n}\in \mathcal{N}_{\var}$, we obtain
$$
M(\|t_{u_{n}}u_{n}\|_{\var}^{2})t_{u_{n}}\|u_{n}\|_{\var}^{2}=\int_{\R^{3}}g(\var x, t_{u_{n}}u_{n})u_{n}, \ \forall \ n\in \N
$$
and passing to the limit $n\rig \infty$, it follows that
$$
M(\|t_{0}u\|_{\var}^{2})t_{0}\|u\|_{\var}^{2}=\int_{\R^{3}}g(\var x, t_{0}u)u,
$$
thence $t_{0}u\in \mathcal{N}_{\var}$ and $t_{u}=t_{0}$, showing that
$\cha{m}_{\var}(u_{n})\rig \cha{m}_{\var}(u)$ in $H_{\var}$. So, $\cha{m}_{\var}$ and $m_{\var}$ are continuous functions and
$(A_{3})$ is proved.

$(A_{4})$ Finally, let $\{u_{n}\}\subset S_{\var}^{+}$ be a sequence such that  $dist(u_{n}, \partial H^{+}_{\var})\rig 0$.
Since, for each $v\in \partial H^{+}_{\var}$ and $n\in \N$, we have
$$
u_{n}^{+}\leq|u_{n}-v| \ \mbox{in} \ \Omega_{\var},
$$
it follows that
\begin{equation}\label{111}
|u_{n}^{+}|_{L^{s}(\Omega_{\var})}\leq\inf_{v\in \partial H^{+}_{\var}}|u_{n}-v|_{L^{s}(\Omega_{\var})}, \ \forall \ n\in \N \ \mbox{and} \ \forall s\in [2, 6].
\end{equation}
Hence, from $(V_{1})$, $(V_{2})$ and Sobolev's embedding, there is
$C(s)>0$ such that
\begin{eqnarray*}
|u_{n}^{+}|_{L^{s}(\Omega_{\var})}&\leq& C(s)\inf_{v\in \partial H^{+}_{\var}}
\left\{\int_{\Omega_{\var}}\left[|\nabla(u_{n}-v)|^{2}+V(\var x)(u_{n}-v)^{2}\right]\right\}^{1/2}\\
&\leq& C(s) dist(u_{n},  \partial H^{+}_{\var}), \ \forall \ n\in \N.
\end{eqnarray*}
From $(g_{1})$, $(g_{2})$ and $(g_{3})(ii)$, there is positive constants $C_{1}$ and $C_{2}$, such that, for each $t>0$
\begin{eqnarray*}
\int_{\R^{3}}G(\var x, tu_{n})&\leq& \int_{\Omega_{\var}}F(tu_{n})+\frac{t^{2}}{K}\int_{\R^{3}\backslash\Omega_{\var}}V(\var x)u_{n}^{2}\\
&\leq& C_{1} t^{4}\int_{\Omega_{\var}}(u_{n}^{+})^{4}+C_{2}t^{q}\int_{\Omega_{\var}}(u_{n}^{+})^{q}+ \frac{1}{K}t^{2}\|u_{n}\|_{\var}^{2}\\
&\leq&  C_{1} C(4)^{4}t^{4} dist(u_{n},  \partial H^{1, +}(\R^{3}))^{4}+C_{2}C(q)t^{q}dist(u_{n},  \partial H^{1, +}(\R^{3}))^{q}+ \frac{1}{K}t^{2}.
\end{eqnarray*}
Therefore,
$$
\limsup_{n\to\infty}\int_{\R^{3}}G(\var x, tu_{n})\leq \frac{1}{K}t^{2}, \ \forall t>0.
$$
From definition of $m_{\var}$, we have
$$
\liminf_{n\to\infty}J_{\var}(m_{\var}(u_{n}))\geq \liminf_{n\to\infty}J_{\var}(t u_{n})\geq \frac{1}{2}\cha{M}(t^{2})-\frac{1}{K}t^{2}, \ \forall \ t>0.
$$
It follows from $(M_{1})$ and from the particular choice of $K$, that
$$
\lim_{n\to\infty}J_{\var}(m_{\infty}(u_{n}))=\infty.
$$
Since $\frac{1}{2}\cha{M}(t_{u_{n}}^{2})\geq J_{\var}(m_{\var}(u_{n}))$, for each $n\in\N$, we conclude from $(M_{3})$ that
$\|m_{\var}(u_{n})\|_{\var}\rig \infty$ as $n\rig \infty$. The Lemma is proved.

$\fim$

\vspace{0.4cm}

We set the applications
$$
\cha{\Psi}_{\var}:H_{\var}^{+}\rig\R \ \mbox{and} \ \Psi_{\var}:S_{\var}^{+}\rig \R,
$$
by $\cha{\Psi}_{\var}(u)=J_{\var}(\cha{m}_{\var}(u))$ \ and $\Psi_{\var}:=(\cha{\Psi}_{\var})_{|_{S_{\var}^{+}}}$.

The next proposition is a direct consequences of the Lemma
\ref{lema3.2.2}. The details can be seen in the relevant material
from $\cite{Szulkin}$. For the convenience of the reader, here we do
 a sketch of the proof.

\begin{Prop}\label{Proposition 4.1.2}
Suppose that the function $M$ satisfies $(M_{1})-(M_{3})$, the potential $V$ satisfies $(V_{1})-(V_{2})$ and the function
$f$ satisfies $(f_{1})-(f_{4})$. Then:
\begin{description}
\item[($a$)] $\cha{\Psi}_{\var}\in C^{1}(H_{\var}^{+}, \R)$ and
$$
\cha{\Psi}_{\var}'(u)v=\frac{\|\cha{m}_{\var}(u)\|_{\var}}{\|u\|_{\var}}J_{\var}'(\cha{m}_{\var}(u))v, \ \forall u\in H_{\var}^{+} \ \mbox{and} \ \forall v\in H_{\var}.
$$

\item[($b$)] $\Psi_{\var}\in C^{1}(S_{\var}^{+}, \R)$ and
$$
\Psi_{\var}'(u)v=\|m_{\var}(u)\|_{\var}J_{\var}'(m_{\var}(u))v, \ \forall v\in T_{u}S_{\var}^{+}.
$$

\item[($c$)] If $\{u_{n}\}$ is a $(PS)_{d}$ sequence by $\Psi_{\var}$ then  $\{m_{\var}(u_{n})\}$ is a $(PS)_{d}$ sequence by $J_{\var}$.
If $\{u_{n}\}\subset \mathcal{N}_{\var}$ is a bounded $(PS)_{d}$ sequence by $J_{\var}$ then $\{m^{-1}_{\var}(u_{n})\}$ is a $(PS)_{d}$
sequence by $\Psi_{\var}$.

\item[($d$)] $u$ is a critical point of $\Psi_{\var}$ if, and only if, $m_{\var}(u)$ is a nontrivial critical point of $J_{\var}$. Moreover, corresponding critical values coincide and
$$
\inf_{S_{\var}^{+}}\Psi_{\var}=\inf_{\mathcal{N}_{\var}}J_{\var}.
$$

\end{description}
\end{Prop}
\noindent {\bf Proof.}  $(a)$ Consider $u\in H_{\var}^{+}$ and $v\in
H_{\var}$. From definition of $\cha{\Psi}_{\var}$, definition of
$t_{u}$ and mean value Theorem,
\begin{eqnarray*}
\cha{\Psi}_{\var}(u+sv)-\cha{\Psi}_{\var}(u)&=&J_{\var}(t_{u+sv}(u+sv))-J_{\var}(t_{u}u)\\
&\leq& J_{\var}(t_{u+sv}(u+sv))-J_{\var}(t_{u+sv}u)\\
&=& J'_{\var}(t_{u+sv}(u+\tau sv))t_{u+sv}sv,
\end{eqnarray*}
where $|s|$ is small sufficient and $\tau\in (0, 1)$. On the other hand,
$$
\cha{\Psi}_{\var}(u+sv)-\cha{\Psi}_{\var}(u)\geq J_{\var}(t_{u}(u+sv))-J_{\var}(t_{u}u)=J'_{\var}(t_{u}(u+\varsigma sv))t_{u}sv,
$$
where $\varsigma\in (0,1)$. Since $u\mapsto t_{u}$ is a continuous application, follows from previous inequalities that
$$
\lim_{s\to 0}\frac{\cha{\Psi}_{\var}(u+sv)-\cha{\Psi}_{\var}(u)}{s}=
t_{u}I'_{\var}(t_{u}u)v=\frac{\|\cha{m}_{\var}(u)\|_{\var}}{\|u\|_{\var}}J'_{\var}(\cha{m}_{\var}(u))v.
$$
Since $J_{\var}\in C^{1}$, it follows that the Gateaux derivative of $\cha{\Psi}_{\var}$ is linear, bounded on $v$
and it is continuous on $u$. From [\cite{Willem}, Proposition 1.3],  $\cha{\Psi}_{\var}\in C^{1}(H_{\var}^{+},\R)$ and
$$
\cha{\Psi}_{\var}'(u)v=\frac{\|\cha{m}_{\var}(u)\|_{\var}}{\|u\|_{\var}}J_{\var}'(\cha{m}_{\var}(u))v,
\ \forall u\in H_{\var}^{+} \ \mbox{and} \ \forall v\in H_{\var}.
$$
The item $(a)$ is proved.

$(b)$ The item $(b)$ is a direct consequences of the item $(a)$.

$(c)$ Once $H_{\var}=T_{u}S_{\var}^{+}\oplus \R \ u$ for each $u\in
S_{\var}^{+}$, the linear projection $P:H_{\var}\rig
T_{u}S_{\var}^{+}$ defined by $P(v+tu)=v$ has uniformly bounded norm
with respect to $u\in S_{\var}^{+}$. Indeed, since $J(u)v:=(u,
v)_{\var}$ is bounded on bounded sets and $J(u)(v+tu)=t$, it follows
that  $\|v+tu\|_{\var}=1$, therefore
$$
|t|\leq \|J(u)\|\|v+tu\|_{\var}=\|J(u)\|\leq C,
$$
for some $C>0$, which implies
$$
\|v\|_{\var}\leq |t|+\|v+tu\|_{\var}\leq (C+1)\|v+tu\|_{\var}, \ \forall \ u\in S_{\var}^{+}, v\in T_{u}S_{\var}^{+} \ \mbox{and} \ t\in \R.
$$
From item $(a)$, we obtain
\begin{equation}\label{17}
\|\Psi'_{\var}(u)\|_{\ast}=\sup_{v\in T_{u}S_{\var}^{+}\atop \|v\|_{\var}=1}\Psi'_{\var}(u)v
=\| w\|_{\var}\sup_{v\in T_{u}S_{\var}^{+}\atop \|v\|_{\var}=1}J'_{\var}(w)v,
\end{equation}
where $w=m_{\var}(u)$. Since $w\in \mathcal{N}_{\var}$, we conclude that $J'_{\var}(w)u=J'_{\var}(w)\frac{w}{\|w\|_{\var}}=0$. By $(a)$
$$
\|\Psi'_{\var}(u)\|_{\ast}\leq \|w\|\|J'_{\var}(w)\|=\|w\|_{\var}\sup_{v\in T_{u}S_{\var}^{+}, t\in \R \atop v+tu\neq 0}
\frac{J'_{\var}(u)(v+tu)}{\|v+tu\|_{\var}}.
$$
Hence,
$$
\|\Psi'_{\var}(u)\|_{\ast}\leq (C+1)\| w\|_{\var}\sup_{v\in
T_{u}S_{\var}^{+}\backslash\{0\}}
\frac{J'_{\var}(w)(v)}{\|v\|_{\var}}=(C+1)\|\Psi'_{\var}(u)\|_{\ast},
$$
showing that,
\begin{equation}\label{16}
\|\Psi'_{\var}(u)\|_{\ast}\leq \|w\|\|J'_{\var}(w)\|\leq (C+1)\|\Psi'_{\var}(u)\|_{\ast}, \ \forall \ u\in S_{\var}^{+}.
\end{equation}
Since $w\in \mathcal{N}{\var}$, we have $\|w\|\geq \tau>0$. Therefore, the inequality in $(\ref{16})$ together
with $J_{\var}(w)=\Psi_{\var}(u)$ imply the item $(c)$.

$(d)$ Follows from $(\ref{17})$ that $\Psi'_{\var}(u)=0$ if, and only if, $J'_{\var}(w)=0$. The remainder follows from
definition of $\Psi_{\var}$.

$\fim$

By using $(M_{1})-(M_{3})$ we have, as in $\cite{Szulkin}$,
the following variational characterization of the infimum of
$J_{\var}$ over $\mathcal{N}_{\var}$:

\begin{eqnarray}\label{minimax}
c_{\var}=\inf_{u\in \mathcal{N}_{\var}}J_{\var}(u)=\inf_{u\in
H^{+}_{\var}(\R^{3})}\max_{t>0}J_{\var}(tu)=\inf_{u\in
S_{\var}^{+}}\max_{t>0}J_{\var}(tu).
\end{eqnarray}

%%%%%%%%%%%%%%%%%%%%%%%%%%%%%%%%%%%%%%%%%%%%%%%%%%%%%%%%%%%%%%%%%%%%%%%%%%%%%%%%%%%%%%%%%%%%%%%%%%%%%%%%%%%%%%%%%%%%%%%%%%%%%%%%%%%%
%%%%%%%%%%%%%%%%%%%%%%%%%%%%%%%%%%%%%%%%%%%%%%%%%%%%%%%%%%%%%%%%%%%%%%%%%%%%%%%%%%%%%%%%%%%%%%%%%%%%%%%%%%%%%%%%%%%%%%%%%%%%%%%%%%%%

The main feature of the modified functional is that it satisfies the
Palais-Smale condition, as we can see from the next results.

\begin{Le}\label{Lemma 4.2.1}
Let $\{u_{n}\}$ be a $(PS)_{d}$ sequence for $J_{\var}$. Then $\{u_{n}\}$ is bounded.
\end{Le}

\noindent {\bf Proof.} Since $\{u_{n}\}$ a $(PS)_{d}$ sequence for
$J_{\var}$, then there is $C>0$ such that
$$
C+\|u_{n}\|_{\var}\geq
J_{\var}(u_{n})-\frac{1}{\theta}J_{\var}'(u_{n})u_{n}, \ \forall
n\in \N.
$$
From $(M_{3})$ and $(g_{3})$, we obtain
$$
C+\|u_{n}\|_{\var}\geq \left(\frac{\theta-2}{2\theta}\right)\left(m_{0}-\frac{1}{K}\right)\|u_{n}\|_{\var}^{2}, \ \forall n\in
\N.
$$
Therefore $\{u_{n}\}$ is bounded in $H_{\var}$.
$\fim$

\begin{Le}\label{Lemma 4.2.2}
Let $\{u_{n}\}$ be a $(PS)_{d}$ sequence for $J_{\var}$. Then for each $\xi>0$, there is $R=R(\xi)>0$ such that
$$
\limsup_{n\to \infty}\int_{\R^{3}\backslash B_{R}}\left[ |\nabla u_{n}|^{2}+V(\var x) u_{n}^{2}\right]<\xi.
$$
\end{Le}

\noindent {\bf Proof.}  Let $\eta_{R} \in C^{\infty}(\R^3)$ such that
\begin{eqnarray*}
\eta_{R}(x) =\ \  \left\{
             \begin{array}{l}
              0 \quad se \quad x \in B_{R/2}(0)
        \\
        \\
        1 \quad se \quad x \not\in B_{R}(0),
        \\
        \\
             \end{array}
           \right.
\end{eqnarray*}
where $0 \leq \eta_{R}(x) \leq 1$, $|\nabla \eta_{R}|
\leq\dis\frac{C}{R}$ and $C$ is a constant independent on $R$. Note
that $\{\eta_{R} u_{n}\}$ is bounded in $H_{\var}$. From definition
of $J_{\var}$
\begin{eqnarray*}
\dis\int_{\R^3}\eta_{R} M(\|u_{n}\|_{\var}^{2})\left[ |\nabla u_n|^{2} + V(\var
x)u_{n}^{2}\right] &=& J'_{\var}(u_n)u_{n}\eta_{R} + \dis\int_{\R^3}g(\var x, u_n)u_{n}\eta_{R}\\
&-&\dis\int_{\R^3}M(\|u_{n}\|_{\var}^{2})u_{n}\nabla u_{n}\nabla \eta_{R}.
\end{eqnarray*}
Choosing $R>0$ such that $\Omega_{\var}\subset B_{\frac{R}{2}}(0)$ and by using $(M_{1})$ and $(g_{3})(ii)$, we have
\begin{eqnarray*}
m_{0}\dis\int_{\R^3}\eta_{R}\left[ |\nabla u_n|^{2} + V(\var
x)u_{n}^{2}
\right] &\leq& J'_{\var}(u_n)u_{n}\eta_{R}\\
&+& \dis\int_{\R^3}\dis\frac{1}{K}V(\var x)u_{n}^{2}\eta_{R} -
\dis\int_{\R^3}M(\|u_{n}\|_{\var}^{2})u_{n}\nabla u_{n}\nabla \eta_{R}.
\end{eqnarray*}
Therefore,
\begin{eqnarray*}
\left(m_{0} -
\dis\frac{1}{K}\right)\dis\int_{\R^3}\eta_{R}\left[|\nabla
u_{n}|^{2} + V(\var x)u_{n}^{2}\right] \leq
|J'_{\var}(u_{n})u_{n}\eta_{R}| + \dis\int_{\R^3}M(\|u_{n}\|_{\var}^{2})u_{n}|\nabla u_{n}\nabla \eta_{R}|.
\end{eqnarray*}
By using Cauchy-Schwarz inequality in $\R^{3}$, Holder's inequality,
the definition of $\eta_{R}$ and from the boundedness of $\{u_{n}\}$
in $H_{\var}$, we conclude that
\begin{eqnarray*}
\dis\int_{\R^3\backslash B_{R}}\left[|\nabla u_n|^{2} +
V(\var x)u_{n}^{2}\right] \leq
C|J'_{\var}(u_{n})u_{n}\eta_{R}| +
\dis\frac{C}{R}\|u_n\|_{\var}.
\end{eqnarray*}
Since $\{u_{n}\}$ and $\{u_{n}\eta_{R}\}$ are bounded in $H_{\var}$, passing to the upper limit of $n\rig \infty$, we obtain
\begin{eqnarray*}
\dis\limsup_{n\rig \infty}\dis\int_{\R^3\backslash
B_{R}}\left[|\nabla u_n|^{2} + V(\var x)u^{2}_{n}\right] \leq
\dis\frac{\widetilde{C}}{R} < \xi,
\end{eqnarray*}
whenever $R=R(\xi)>\frac{\widetilde{C}}{\xi}$.
$\fim$

The next result does not appear in \cite{Del Pino}, however, since that we are
working with the Kirchhoff problem type, it is required here.

\begin{Le}\label{Lemma 4.2.3}
Let $\{u_{n}\}$ be a $(PS)_{d}$ sequence for $J_{\var}$ such that
$u_{n}\righ u$, then
$$
\lim_{n\to \infty}\int_{B_{R}}\left[ |\nabla u_{n}|^{2}+V(\var x)
u_{n}^{2}\right]=\int_{B_{R}}\left[ |\nabla u|^{2}+V(\var x)
u^{2}\right],
$$
for all $R>0$.
\end{Le}

\noindent {\bf Proof.} We can assume that $\|u_{n}\|_{\var}\rig t_{0}$, thus, we have $\|u\|_{\var}\leq t_{0}$. Let $\eta_{\rho} \in C^{\infty}(\R^3)$ such that
\begin{eqnarray*}
\eta_{\rho}(x) =\ \  \left\{
             \begin{array}{l}
              1 \quad se \quad x \in B_{\rho}(0)
        \\
        \\
        0 \quad se \quad x \not\in B_{2\rho}(0).
        \\
        \\
             \end{array}
           \right.
\end{eqnarray*}
with $0 \leq \eta_{\rho}(x) \leq 1$. Let,
$$
P_{n}(x)=M(\|u_{n}\|_{\var}^{2})\left[ |\nabla u_{n}-\nabla u|^{2}+V(\var x)(u_{n}-u)^{2}\right].
$$
For each $R>0$ fixed, choosing $\rho>R$ we obtain
$$
\int_{B_{R}}P_{n}=\int_{B_{R}}P_{n}\eta_{\rho}\leq
M(\|u_{n}\|_{\var}^{2})\int_{\R^{3}}\left[ |\nabla u_{n}-\nabla
u|^{2}+V(\var x)(u_{n}-u)^{2}\right]\eta_{\rho}.
$$
By expanding the inner product in $\R^{3}$,
\begin{eqnarray*}
\int_{B_{R}}P_{n}&\leq& M(\|u_{n}\|_{\var}^{2})\int_{\R^{3}}\left[ |\nabla u_{n}|^{2}+V(\var x)(u_{n})^{2}\right]\eta_{\rho}\\
&-&2M(\|u_{n}\|_{\var}^{2})\int_{\R^{3}}\left[ \nabla u_{n}\nabla u+V(\var x)u_{n}u\right]\eta_{\rho}\\
&+&M(\|u_{n}\|_{\var}^{2})\int_{\R^{3}}\left[ |\nabla u|^{2}+V(\var x)u^{2}\right]\eta_{\rho}.
\end{eqnarray*}
Setting
$$
I_{n,\rho}^{1}=M(\|u_{n}\|_{\var}^{2})\int_{\R^{3}}\left[ |\nabla u_{n}|^{2}+V(\var x)(u_{n})^{2}\right]\eta_{\rho}
-\int_{\R^{3}}g(\var x, u_{n})u_{n}\eta_{\rho},
$$
$$
I_{n,\rho}^{2}=M(\|u_{n}\|_{\var}^{2})\int_{\R^{3}}\left[ \nabla u_{n}\nabla u+V(\var x)u_{n}u\right]\eta_{\rho}
-\int_{\R^{3}}g(\var x, u_{n})u\eta_{\rho},
$$
$$
I_{n,\rho}^{3}=-M(\|u_{n}\|_{\var}^{2})\int_{\R^{3}}\left[ \nabla u_{n}\nabla u+V(\var x)u_{n}u\right]\eta_{\rho}
+M(\|u_{n}\|_{\var}^{2})\int_{\R^{3}}\left[ |\nabla u|^{2}+V(\var x)u^{2}\right]\eta_{\rho}
$$
and
$$
I_{n,\rho}^{4}=\int_{\R^{3}}g(\var x, u_{n})u_{n}\eta_{\rho}-\int_{\R^{3}}g(\var x, u_{n})u\eta_{\rho}.
$$
We have that,
\begin{equation}\label{11}
0\leq\int_{B_{R}}P_{n}\leq |I_{n,\rho}^{1}|+|I_{n,\rho}^{2}|+|I_{n,\rho}^{3}|+|I_{n,\rho}^{4}|.
\end{equation}
Observe that
$$
I_{n,\rho}^{1}=J_{\var}'(u_{n})(u_{n}\eta_{\rho})-M(\|u_{n}\|_{\var}^{2})\int_{\R^{3}} u_{n}\nabla u_{n}\nabla\eta_{\rho}.
$$
Since $\{u_{n}\eta_{\rho}\}$ is bounded in $H_{\var}$, we have
$J_{\var}'(u_{n})(u_{n}\eta_{\rho})=o_{n}(1)$. Moreover, from a
straightforward computation
$$
\lim_{\rho\to \infty}\left[ \limsup_{n\to \infty}
\left|M(\|u_{n}\|_{\var}^{2})\int_{\R^{3}} u_{n}\nabla
u_{n}\nabla\eta_{\rho}\right|\right]=0.
$$
Then,
\begin{equation}\label{12}
\lim_{\rho\to \infty}\left[ \limsup_{n\to \infty} |I_{n,\rho}^{1}|\right]=0.
\end{equation}

We see also that
$$
I_{n,\rho}^{2}=J_{\var}'(u_{n})(u\eta_{\rho})-M(\|u_{n}\|_{\var}^{2})\int_{\R^{3}} u\nabla u_{n}\nabla\eta_{\rho}.
$$
By arguing in the same way as in the previous case,
$$
J_{\var}'(u_{n})(u\eta_{\rho})=o_{n}(1)
$$
and
$$
\lim_{\rho\to \infty}\left[ \limsup_{n\to \infty} \left|M(\|u_{n}\|_{\var}^{2})\int_{\R^{3}} u\nabla u_{n}\nabla\eta_{\rho}\right|\right]=0.
$$
Therefore,
\begin{equation}\label{13}
\lim_{\rho\to \infty}\left[ \limsup_{n\to \infty} |I_{n,\rho}^{2}|\right]=0.
\end{equation}

On the other hand, from the weak convergence
\begin{equation}\label{14}
\lim_{n\to \infty} |I_{n,\rho}^{3}|=0, \ \forall \ \rho>R.
\end{equation}

Finally, from
$$
u_{n}\rig u, \ \mbox{in} \ L^{s}_{loc}(\R^{3}), 1\leq s<6,
$$
we conclude that
\begin{equation}\label{15}
\lim_{n\to \infty} |I_{n,\rho}^{4}|=0, \ \forall \ \rho>R.
\end{equation}
From $(\ref{11})$, $(\ref{12})$, $(\ref{13})$, $(\ref{14})$ and $(\ref{15})$, we obtain
$$
0\leq \limsup_{n\to \infty} \int_{B_{R}}P_{n}\leq 0.
$$
Hence, $\dis\lim_{n\to \infty} \dis\int_{B_{R}}P_{n}=0$ and consequently
$$
\lim_{n\to \infty}\int_{B_{R}}\left[ |\nabla u_{n}|^{2}+V(\var x) u_{n}^{2}\right]=\int_{B_{R}}\left[ |\nabla u|^{2}+V(\var x) u^{2}\right].
$$
$\fim$

\begin{Prop}\label{Proposition 4.2.1}
The functional $J_{\var}$ verifies the $(PS)_{d}$ condition in $H_{\var}$.
\end{Prop}

\noindent {\bf Proof.} Let $\{u_{n}\}$ be a $(PS)_{d}$ sequence for $J_{\var}$. From Lemma $\ref{Lemma 4.2.1}$ we know that $\{u_{n}\}$ is bounded in $H_{\var}$.
Passing to a subsequence, we obtain
$$
u_{n}\righ u, \ \mbox{in} \ H_{\var}.
$$
From Lemma $\ref{Lemma 4.2.2}$, it follows that for each $\xi>0$ given there is $R=R(\xi)>\frac{\widetilde{C}}{\xi}$ with $\widetilde{C}$ independent on $\xi$ such that
$$
\limsup_{n\to \infty}\int_{\R^{3}\backslash B_{R}}\left[ |\nabla u_{n}|^{2}+V(\var x) u_{n}^{2}\right]<\xi.
$$
Therefore, from Lemma $\ref{Lemma 4.2.3}$,

\begin{eqnarray*}
\|u\|_{\var}^{2}&\leq&\liminf_{n\to \infty}\|u_{n}\|_{\var}^{2}
\leq\limsup_{n\to \infty}\|u_{n}\|_{\var}^{2}\\
&=&\limsup_{n\to \infty}\left\{\int_{B_{R}}\left[ |\nabla u_{n}|^{2}+V(\var x) u_{n}^{2}\right]+\int_{\R^{3}\backslash B_{R}}\left[ |\nabla u_{n}|^{2}+V(\var x) u_{n}^{2}\right]\right\}\\
&=&\int_{B_{R}}\left[ |\nabla u|^{2}+V(\var x) u^{2}\right]+\limsup_{n\to \infty}\int_{\R^{3}\backslash B_{R}}\left[ |\nabla u_{n}|^{2}+V(\var x) u_{n}^{2}\right]\\
&<&\int_{B_{R}}\left[ |\nabla u|^{2}+V(\var x) u^{2}\right]+\xi,
\end{eqnarray*}
where $R=R(\xi)>\frac{\widetilde{C}}{\xi}$. Passing to the limit of
$\xi\to 0$ we have $R\rig\infty$, which implies
$$
\|u\|_{\var}^{2}\leq \liminf_{n\to \infty}\|u_{n}\|_{\var}^{2}\leq\limsup_{n\to \infty}\|u_{n}\|_{\var}^{2}\leq \|u\|_{\var}^{2},
$$
and so $\|u_{n}\|_{\var}\rig \|u\|_{\var}$ and consequently $u_{n}\rig u$ in $H_{\var}$.
$\fim$

Since $f$ is only continuous and $V$ has geometry of the Del Pino
and Felmer type \cite{Del Pino}, in the next result (which is required for the multiplicity result)
we use arguments that don't appear in \cite{Del Pino} and \cite{Wang 1}.

\begin{Co}\label{Corollary 4.2.1}
The functional $\Psi_{\var}$ verifies the $(PS)_{d}$ condition on $S_{\var}^{+}$.
\end{Co}

\noindent {\bf Proof.} Let $\{u_{n}\}\subset S_{\var}^{+}$ be a $(PS)_{d}$ sequence for $\Psi_{\var}$. Thus,
$$
\Psi_{\var}(u_{n})\rig d
$$
and
$$
\left\|\Psi_{\var}'(u_{n})\right\|_{\ast}\rig 0,
$$
where $\|.\|_{\ast}$ is the norm in the dual space $\left(T_{u_{n}}S_{\var}^{+}\right)'$. From
[$\cite{Willem}$, Proposition 5.12] there is a sequence $\{\lambda_{n}\}\subset \R$ such that
$$
\cha{\Psi}_{\var}'(u_{n})v=2\lambda_{n}(u_{n}, v)_{\var}+o_{n}(1), \
\forall v\in H_{\var}.
$$
From item $(a)$ of the Proposition $\ref{Proposition 4.1.2}$,
$$
\|m_{\var}(u_{n})\|_{\var}J_{\var}'(m_{\var}(u_{n}))v=2\lambda_{n}(u_{n}, v)_{\var}+o_{n}(1), \ \forall v\in H_{\var}.
$$
So,
\begin{equation}
J_{\var}'(m_{\var}(u_{n}))v=2\lambda_{n}\left(u_{n}, \frac{v}{\|m_{\var}(u_{n})\|_{\var}}\right)_{\var}+o_{n}(1), \ \forall v\in H_{\var}.\label{equation 4.2.1}
\end{equation}
Putting $v=m_{\var}(u_{n})$, we have
$$
0=J_{\var}'(m_{\var}(u_{n}))m_{\var}(u_{n})=2\lambda_{n}+o_{n}(1),
$$
because $\{u_{n}\}\subset S_{\var}^{+}$ yields $u_{n}=\frac{m_{\var}(u_{n})}{\|m_{\var}(u_{n})\|_{\var}}$.
It shows that $\lambda_{n}\rig 0$ as $n\rig\infty$. From $(\ref{equation 4.2.1})$,
$$
|J_{\var}'(m_{\var}(u_{n}))v|\leq 2\frac{\|v\|_{\var}}{\|m_{\var}(u_{n})\|_{\var}}|\lambda_{n}|+o_{n}(1),
$$
and from Lemma $\ref{lema3.2.2}(A_{2})$,
$$
\|J_{\var}'(m_{\var}(u_{n}))\|\leq C|\lambda_{n}|+o_{n}(1)=o_{n}(1).
$$
Therefore $\{m_{\var}(u_{n})\}$ is a $(PS)_{d}$ sequence for $J_{\var}$ in $H_{\var}$ and from Proposition $\ref{Proposition 4.2.1}$ we obtain $\overline{u}\in H_{\var}$ such that $m_{\var}(u_{n})\rig \overline{u}$, with $\|\overline{u}\|_{\var}\geq \tau>0$. Hence,
$$
u_{n}=\frac{m_{\var}(u_{n})}{\|m_{\var}(u_{n})\|_{\var}}\rig \frac{\overline{u}}{\|\overline{u}\|_{\var}}, \ \mbox{in} \ H_{\var}.
$$
$\fim$

\begin{Th}\label{Theorem 4.2.1}
Suppose that the function $M$ satisfies $(M_{1})-(M_{3})$, the
potential $V$ satisfies $(V_{1})-(V_{2})$ and the function $f$ satisfies
$(f_{1})-(f_{4})$. Then, the auxiliary problem $(P_{\var,A})$ has a
positive ground-state solution for all $\var>0$.
\end{Th}

\noindent {\bf Proof.} This result follows from Lemma \ref{Lemma
4.1.1},  Proposition \ref{Proposition 4.2.1} and maximum
principle.$\fim$

%%%%%%%%%%%%%%%%%%%%%%%%%%%%%%%%%%%%%%%%%%%%%%%%%%%%%%%%%%%%%%%%%%%%%%%%%%%%%%%%%%%%%%%%%%%%%%%%%%%%%%%%%%%%%%%%%%%%%%%%%%%%%%%%%%%%%%%%%%
%%%%%%%%%%%%%%%%%%%%%%%%%%%%%%%%%%%%%%%%%%%%%%%%%%%%%%%%%%%%%%%%%%%%%%%%%%%%%%%%%%%%%%%%%%%%%%%%%%%%%%%%%%%%%%%%%%%%%%%%%%%%%%%%%%%%%%%%%%
%%%%%%%%%%%%%%%%%%%%%%%%%%%%%%%%%%%%%%%%%%%%%%%%%%%%%%%%%%%%%%%%%%%%%%%%%%%%%%%%%%%%%%%%%%%%%%%%%%%%%%%%%%%%%%%%%%%%%%%%%%%%%%%%%%%%%%%%%%

\section{Multiplicity of solutions of auxiliary problem}

\subsection{The autonomous problem}

Since we are interested in giving a multiplicity result for the
auxiliary problem, we start by considering the limit problem
associated to $(\widetilde{P}_{\var})$, namely, the problem

$$
 \left\{
\begin{array}{rcl}
\mathfrak{L}_{0}u=f(u), \ \R^{3}\\
u>0, \ \R^{3}\\
u \in H^{1}(\R^3)
\end{array}
\right.\leqno{(P_{0})}
$$
where

$$
\mathfrak{L}_{0}u=M\left(\dis\int_{\R^{3}}|\nabla
u|^{2}+\dis\int_{\R^{3}}V_{0} u^{2}\right)\left[-\Delta u + V_{0} u
\right],
$$
which has the following associated functional
$$
I_{0}(u)=\frac{1}{2}\cha{M}\left(\dis\int_{\R^{3}}|\nabla
u|^{2}+\dis\int_{\R^{3}}V_{0} u^{2}\right)-\dis\int_{\R^{3}}F(u).
$$
This functional is well defined on the Hilbert space
$H_{0}=H^{1}(\R^{3})$ with the inner product
$$
(u,v)_{0}=\dis\int_{\R^{3}}\nabla
u\nabla v+\dis\int_{\R^{3}}V_{0} u v
$$
and norm
$$
\| u\|_{0}^{2}=\dis\int_{\R^{3}}|\nabla
u|^{2}+\dis\int_{\R^{3}}V_{0} u^{2}
$$
fixed. We denote the Nehari manifold associated to $I_{0}$ by
$$
\mathcal{N}_{0}=\{u\in H_{0}\backslash\{0\}: I'_{0}(u)u=0\}.
$$
We denote by $H_{0}^{+}$ the open subset of $H_{0}$
given by
$$
H_{0}^{+}=\{u\in H_{0}: |supp (u^{+})|>0 \},
$$
and $S_{0}^{+}=S_{0}\cap H_{0}^{+}$, where
$S_{0}$ is the unit sphere of $H_{0}$.

As in the section 2, $S_{0}^{+}$ is a incomplete
$C^{1,1}$-manifold of codimension $1$, modeled on $H_{0}$ and contained  in the open
$H_{0}^{+}$. Thus, $H_{0}=T_{u}S_{0}^{+}\oplus
\R \ u$ for each $u\in S_{0}^{+}$, where
$T_{u}S_{0}^{+}=\{v\in H_{0}:(u, v)_{0}=0\}$.

In the sequel we enunciate without proof one Lemma and one
Proposition, which allow us to prove the Lemma $\ref{compacity Lemma}$. The proofs follow from a similar argument to that used in
the proofs of Lemma \ref{lema3.2.2} and Proposition \ref{Proposition 4.1.2}.
\begin{Le}\label{Lemma 3.1.1}
Suppose that the function $M$ satisfies $(M_{1})-(M_{3})$ and the function
$f$ satisfies $(f_{1})-(f_{4})$. So:
\begin{description}

\item[($A_{1}$)] For each $u\in H_{0}^{+}$, let
$h:\R_{+}\rig\R$ be defined by $h_{u}(t)=I_{0}(tu)$. Then,
there is a unique $t_{u}>0$ such that $h_{u}'(t)>0$ in $(0,t_{u})$
and $h_{u}'(t)<0$ in $(t_{u}, \infty)$.

\item[($A_{2}$)] there is $\tau>0$ independent on $u$ such that $t_{u}\geq \tau$ for all $u\in S^{+}_{0}$. Moreover,
for each compact set $\mathcal{W}\subset S^{+}_{0}$ there is
$C_{\mathcal{W}}>0$ such that $t_{u}\leq C_{\mathcal{W}}$, for all
$u\in \mathcal{W}$.

\item[($A_{3}$)] The map
$\widehat{m}:H_{0}^{+}\rightarrow
\mathcal{N}_{0}$ given by $\widehat{m}(u)=t_{u}u$ is
continuous and
$m:=\widehat{m}_{\bigl|S^{+}_{0}}$ is a
homeomorphism between $S^{+}_{0}$ and $\mathcal{N}_{0}$.
Moreover, $m^{-1}(u)= \frac{u}{\|u\|_{0}}$.

\item[($A_{4}$)] If there is a sequence $(u_{n})\subset
S^{+}_{0}$ such that \mbox{dist}$(u_{n},\partial H^{+}_{0})\rightarrow 0$, then
$\|m(u_{n})\|_{0}\rightarrow \infty$ and
$I_{0}(m(u_{n}))\rightarrow \infty$.
\end{description}

\end{Le}

We set the applications
$$
\cha{\Psi}_{0}:H_{0}^{+}\rig\R \ \mbox{and} \
\Psi_{0}:S_{0}^{+}\rig \R,
$$
by $\cha{\Psi}_{0}(u)=I_{0}(\cha{m}(u))$ \ and
$\Psi_{0}:=(\cha{\Psi}_{0})_{|_{S_{0}^{+}}}$.

\begin{Prop}\label{Proposition 3.1.2}
Suppose that the function $M$ satisfies $(M_{1})-(M_{3})$ and the function
$f$ satisfies $(f_{1})-(f_{4})$. So:
\begin{description}
\item[($a$)] $\cha{\Psi}_{0}\in C^{1}(H_{0}^{+}, \R)$ and
$$
\cha{\Psi}_{0}'(u)v=\frac{\|\cha{m}(u)\|_{0}}{\|u\|_{0}}I_{0}'(\cha{m}(u))v,
\ \forall u\in H_{0}^{+} \ \mbox{and} \ \forall v\in
H_{0}.
$$

\item[($b$)] $\Psi_{0}\in C^{1}(S_{0}^{+}, \R)$ and
$$
\Psi_{0}'(u)v=\|m(u)\|_{0}I_{0}'(m(u))v, \ \forall v\in
T_{u}S_{0}^{+}.
$$

\item[($c$)] If $\{u_{n}\}$ is a $(PS)_{d}$ sequence by $\Psi_{0}$ then  $\{m(u_{n})\}$ is a $(PS)_{d}$ sequence by $I_{0}$.
If $\{u_{n}\}\subset \mathcal{N}_{0}$ is a bounded $(PS)_{d}$ sequence by $I_{0}$ then $\{m^{-1}(u_{n})\}$ is a $(PS)_{d}$ sequence by $\Psi_{0}$.

\item[($d$)] $u$ is a critical point of $\Psi_{0}$ if, and only if, $m(u)$ is a nontrivial critical point of $I_{0}$.
Moreover, corresponding critical values coincide and
$$
\inf_{S_{0}^{+}}\Psi_{0}=\inf_{\mathcal{N}_{0}}I_{0}.
$$

\end{description}
\end{Prop}

\begin{Remark}
As in the section 2, there holds
\begin{eqnarray}\label{charactetization minimax}
c_{0}=\inf_{u\in \mathcal{N}_{0}}I_{0}(u)=\inf_{u\in
H^{+}_{0}}\max_{t>0}I_{0}(tu)=\inf_{u\in
S_{0}^{+}}\max_{t>0}I_{0}(tu).
\end{eqnarray}
\end{Remark}

The next Lemma allows us to assume that the weak limit of a
$(PS)_{d}$ sequence is non-trivial.

\begin{Le}\label{Lemma 3.3}
Let $\{u_{n}\}\subset H_{0}$ be a $(PS)_{d}$ sequence for $I_{0}$ with
$u_{n}\righ 0$. Then, only one of the alternatives below hold:

a) $u_{n}\rig 0$ in $H_{0}$

b) there is a sequence $({y_{n}})\subset \R^{3}$ and constants $R, \beta>0$
such that
$$
\liminf_{n\to \infty}\int_{B_{R}(y_{n})}u_{n}^{2}\geq \beta>0.
$$

\end{Le}

\noindent {\bf Proof.} Suppose that b) doesn't hold. It follows that
for all $R>0$ we have
$$
\lim_{n\to \infty}\sup_{y\in \R^{3}}\int_{B_{R}(y)}u_{n}^{2}=0.
$$
Since $\{u_{n}\}$ is bounded in $H_{0}$, we conclude from
[$\cite{Willem}$, Lemma 1.21] that

$$
u_{n}\rig 0 \ \mbox{in} \ L^{s}(\R^{3}), 2<s<6.
$$
From $(M_{1})$, $(f_{1})$ and $(f_{2})$,
$$
0\leq m_{0}\|u_{n}\|_{0}\leq\dis\int_{\R^{3}}
f(u_{n})u_{n}+o_{n}(1)=o_{n}(1).
$$
Therefore the item $a)$ is true. $\fim$
\begin{Remark}\label{limitefraconaonulo}
As it has been mentioned, if $u$ is the weak limit of a $(PS)_{c_{0}}$
sequence $\{u_{n}\}$ for the functional $I_{0}$, then we can assume
$u\neq 0$, otherwise we would have $u_{n}\righ 0$ and, once it
doesn't occur $u_{n}\rig 0$, we conclude from the Lemma $\ref{Lemma
3.3}$ that there are $\{y_{n}\}\subset \R^{3}$ and $R,\beta>0$ such
that
$$
\liminf_{n\to \infty}\int_{B_{R}(y_{n})}u_{n}^{2}\geq \beta>0.
$$
Set $v_{n}(x)=u_{n}(x+y_{n})$, making a change of variable, we can
prove that $\{v_{n}\}$ is a $(PS)_{c_{0}}$ sequence for the
functional $I_{0}$, it is bounded in $H_{0}$ and there is $v\in
H_{0}$ with $v_{n}\rightharpoonup v$ in $H_{0}$ with $v\neq 0$.
\end{Remark}

In the next Proposition we obtain a positive ground-state solution
for the autonomous problem $(P_{0})$.

\begin{Th}\label{Theorem 3.1}
Let $\{u_{n}\}\subset H_{0}$ be a $(PS)_{c_{0}}$ sequence for
$I_{0}$. Then there is $u\in H_{0}\backslash\{0\}$ with $u\geq 0$
such that, passing a subsequence, we have $u_{n}\rig u$ in $H_{0}$.
Moreover, $u$ is a positive ground-state solution for the problem
$(P_{0})$.

\end{Th}

\noindent {\bf Proof.} Arguing as Lemma \ref{Lemma 4.2.1}, we have
that $\{u_{n}\}$ is bounded in $H_{0}$. Thus, passing a subsequence
if necessary, we obtain
\begin{equation}
u_{n}\righ u \ \mbox{em} \ H_{0},\label{weak convergence 1}
\end{equation}
\begin{equation}
u_{n}\rig u \ \mbox{em} \ L^{s}_{loc}(\R^{3}), 1\leq
s<6\label{Lebesgue convergence 1}
\end{equation}
and
\begin{equation}
\|u_{n}\|_{0}\rig t_{0}.\label{norm convergence 1}
\end{equation}
So, from $(\ref{weak convergence 1})$ we conclude that
\begin{equation}
(u_{n},v)_{0}\rig (u, v)_{0}, \ \forall v\in H_{0}.\label{inner
product}
\end{equation}
On the other hand, due to density of $C_{0}^{\infty}(\R^{3})$ in
$H_{0}$ and from convergence in $(\ref{Lebesgue convergence 1})$,
it results that
\begin{equation}
\dis\int_{\R^{3}} f(u_{n})v\rig\dis\int_{\R^{3}} f(u)v, \forall v\in
H_{0}.\label{function convergence}
\end{equation}
Now, from convergence in $(\ref{weak convergence 1})$ and
$(\ref{norm convergence 1})$, occurs
$$
\|u\|_{0}^{2}\leq \liminf_{n\to\infty}\|u_{n}\|_{0}^{2}=t_{0}^{2},
$$
and from $(M_{2})$ it follows that $M(\|u\|_{0}^{2})\leq
M(t_{0}^{2})$.

Since $(M_{3})$ implies that the function
$t\mapsto\frac{1}{2}\cha{M}(t)-\frac{1}{4}M(t)t$ is non-decreasing,
we can  argue as in $\cite{Alves and Figueiredo}$ and to prove that
$M(t_{0}^{2})=M(\|u\|_{0}^{2})$ and the theorem now follows. $\fim$

\vspace*{.4cm}

\begin{Remark}\label{Remark 3}
Since functional $I_{0}$ has the mountain pass geometry, it follows from \mbox{[\cite{Willem}, Theorem 1.15]} and Theorem $\ref{Theorem 3.1}$
that $(P_{0})$ admits a positive ground-state solution.
\end{Remark}

The next lemma is a compactness result on the autonomous problem
which we will use later.

\begin{Le}
Let $\{u_{n}\}$ be a sequence in $H^{1}(\R^{3})$ such that
$I_{0}(u_{n})\rig c_{0}$ and $\{u_{n}\} \subset \mathcal{N}_{0}$.
Then, $\{u_{n}\}$ has a convergent subsequence in $H^{1}(\R^{3})$.
\label{compacity Lemma}
\end{Le}

\noindent {\bf Proof.} Since $\{u_{n}\}\subset \mathcal{N}_{0}$, it
follows from item $(A_{3})$ of the Lemma $\ref{Lemma 3.1.1}$,
from item $(d)$ of the Proposition $\ref{Proposition 3.1.2}$ and from
the Remark $\ref{charactetization minimax}$ that
\begin{equation}
v_{n}=m^{-1}(u_{n})=\frac{u_{n}}{\|u_{n}\|_{0}}\in S_{0}^{+}, \
\forall n\in \N\label{equality 1}
\end{equation}
and
$$
\Psi_{0}(v_{n})=I_{0}(u_{n})\rig c_{0}=\inf_{S_{0}^{+}}\Psi_{0}.
$$
Although $S_{0}^{+}$ is incomplete, due to item $(A_{4})$
from the Lemma $\ref{Lemma 3.1.1}$, we can still apply the Ekeland's
variational principle [$\cite{Ekeland}$, Theorem 1.1] to the functional $\xi_{0}: V\rightarrow \R\cup\{\infty\}$
defined by $\xi_{0}(u)=\cha{\Psi}_{0}(u)$ if $u\in H_{0}^{+}$ and $\xi_{0}(u)=\infty$ if $u\in\partial H_{0}^{+}$,
where $V=\overline{H_{0}^{+}}$ is a complete metric space equipped with the metric $d(u,v)=\|u-v\|_{0}$.
In fact, from Lemma $\ref{Lemma 3.1.1} (A_{4})$, $\xi_{0} \in C(V,\R\cup\{\infty\})$ and, from Proposition
$\ref{Proposition 3.1.2}(d)$, $\xi_{0}$  is bounded below. Thus, we can conclude there is a
sequence $\{\widehat{v}_{n}\}\subset S_{0}^{+}$ such that
$\{\widehat{v}_{n}\}$ is a $(PS)_{c_{0}}$ sequence for $\Psi_{0}$ on
$S_{0}^{+}$ and
\begin{eqnarray}\label{proxima}
\|\widehat{v}_{n}-v_{n}\|_{0}=o_{n}(1).
\end{eqnarray}
The remainder of the proof follows by using Proposition $\ref{Proposition 3.1.2}$,
Theorem $\ref{Theorem 3.1}$ and
arguing as in the proof of Corollary $\ref{Corollary 4.2.1}$.

$\fim$

In this section we will relate the number of positive solutions of $(P_{\var,A})$ to topology of $\Pi$, for this we need some preliminary results.

\subsection{Technical results}

Let $\delta>0$ fixed and $\Pi_{\delta} \subset \Omega$. Let $\eta
\in C_{0}^{\infty}([0,\infty))$ be such that $0 \leq \eta(t) \leq
1$, $\eta(t) = 1$ if $0 \leq t \leq \delta/2 $ and $\eta(t) = 0$ if
$t \geq \delta$. We denote by $w$ a positive ground-state solution
of the problem $(P_{0})$ (see Remark $\ref{Remark 3}$).

For each $y \in \Pi= \{x \in \Omega : V(x) = V_{0}\}$, we define the
function
\begin{eqnarray*}
\widetilde{\Upsilon}_{\var,y}(x) = \eta(|\var x -
y|)w\left(\frac{\var x - y}{\var}\right).
\end{eqnarray*}

\noindent Let $t_{\var}> 0$ be the unique positive number such that
\begin{eqnarray*}
\dis\max_{t \geq 0}J_{\var}(t\widetilde{\Upsilon}_{\var,y}) =
J_{\var}(t_{\var}\widetilde{\Upsilon}_{\var,y}).
\end{eqnarray*}

\noindent By noticing that $t_{\var}\widetilde{\Upsilon}_{\var,y}
\in \mathcal{N}_{\var}$, we can  now define the continuous function
\begin{eqnarray*}
\Upsilon_{\var} : & \Pi& \longrightarrow \mathcal{N}_{\var}  \\
          &y & \longmapsto \Upsilon_{\var}(y)=t_{\var}\widetilde{\Upsilon}_{\var,y}.
\end{eqnarray*}

\begin{Le}\label{Lemma 5.0.1}
Let $\Pi \subset \Omega$. Then,
\begin{eqnarray*}
\dis\lim_{\var \rig 0}J_{\var}(\Upsilon_{\var}(y)) = c_{0} \
\mbox{uniformly in} \ y \in \Pi.
\end{eqnarray*}
\end{Le}

\noindent {\bf Proof.}  Arguing by contradiction, we suppose that
there exist  $\delta_{0} > 0$ and a sequence $\{y_n\} \subset \Pi$
verifying
\begin{eqnarray}\label{equation 2}
\mid J_{\var_{n}}(\Upsilon_{\var_{n}}(y_{n})) - c_{0} \mid \geq
\delta _{0}\ \mbox{where} \ \var_{n}\rig 0 \ \mbox{when} \ \ n\rig
\infty.
\end{eqnarray}
\noindent From definition of $\Upsilon_{\var_{n}}(y_{n})$, we have
\begin{equation}
J_{\var_{n}}(\Upsilon_{\var_{n}}(y_{n}))=\dis\frac{1}{2}\cha{M}\left(t^{2}_{\var_{n}}\|
\widetilde{\Upsilon}_{\var_{n},y_{n}}\|_{\var_{n}}^{2}\right)
-\dis\int_{\R^3}G\left(\var_{n}x,t_{\var_{n}}\widetilde{\Upsilon}_{\var_{n},
y_{n}}\right)\label{functional}
\end{equation}
and
\begin{equation}
J_{\var_{n}}'(\Upsilon_{\var_{n}}(y_{n}))\Upsilon_{\var_{n}}(y_{n})=0.\label{equation 3}
\end{equation}
Using definition of $\Upsilon_{\var_{n}}(y_{n})$ again and making the change of variable $z=\frac{\var_{n}x-y_{n}}{\var_{n}}$, we have
\begin{eqnarray*}
&&J_{\var_{n}}(\Upsilon_{\var_{n}}(y_{n}))=\\
&&\dis\frac{1}{2}\cha{M}\left(t^{2}_{\var_{n}}\left(\int_{\R^{3}}\left|\nabla\left(\eta(|\var_{n}z|)w(z)\right)\right|^{2}+
\int_{\R^{3}}V(\var_{n}z+y_{n})\left(\eta(|\var_{n}z|)w(z)\right)^{2}\right)\right)\\ &&-\dis\int_{\R^3}G\left(\var_{n}z+y_{n},t_{\var_{n}}\eta(|\var_{n}z|)w(z)\right).
\end{eqnarray*}
Moreover, putting
$$\Lambda_{n}^{2}=\int_{\R^{3}}\left|\nabla\left(\eta(|\var_{n}z|)w(z)\right)\right|^{2}+
\dis\int_{\R^{3}}V(\var_{n}z+y_{n})\left(\eta(|\var_{n}z|)w(z)\right)^{2},
$$
the equality in $(\ref{equation 3})$ yields
$$
\frac{M(t^{2}_{\var_{n}}\Lambda_{n}^{2})}{t^{2}_{\var_{n}}\Lambda_{n}^{2}}=\frac{1}{\Lambda_{n}^{4}}\dis\int_{\R^3} \left[\frac{g(\var_{n}z+y_{n},t_{\var_{n}}\eta(|\var_{n}z|)w(z))}{\left(t_{\var_{n}}\eta(|\var_{n}z|)w(z)\right)^{3}}\right]
(\eta(|\var_{n}z|)w(z))^{4}.
$$
For each $n\in \N$ and for all $z\in B_{\frac{\delta}{\var_{n}}}(0)$, we have $\var_{n}z\in  B_{\delta}(0)$. So,
$$
\var_{n}z+y_{n}\in B_{\delta}(y_{n})\subset \Pi_{\delta}\subset \Omega.
$$
Since $G=F$ in $\Omega$, it follows from $(\ref{functional})$ that
\begin{equation}
J_{\var_{n}}(\Upsilon_{\var_{n}}(y_{n}))=
\dis\frac{1}{2}\cha{M}(t^{2}_{\var_{n}}\Lambda_{n}^{2})-\dis\int_{\R^3}F\left(t_{\var_{n}}\eta(|\var_{n}z|)w(z)\right)\label{equation
7}
\end{equation}
and
\begin{equation}
\frac{M(t^{2}_{\var_{n}}\Lambda_{n}^{2})}{t^{2}_{\var_{n}}\Lambda_{n}^{2}}=\frac{1}{\Lambda_{n}^{4}}\dis\int_{\R^3} \left[\frac{f(t_{\var_{n}}\eta(|\var_{n}z|)w(z))}{\left(t_{\var_{n}}\eta(|\var_{n}z|)w(z)\right)^{3}}\right]
(\eta(|\var_{n}z|)w(z))^{4}.\label{equation 5}
\end{equation}
From the Lebesgue's theorem, when $n\rig \infty$
\begin{eqnarray}
\|\widetilde{\Upsilon}_{\var_{n},y_{n}}\|_{\var_{n}}^{2}=\Lambda_{n}^{2}\rig
\|w\|_{0}^{2}, \label{equation 4}
\end{eqnarray}
$$
\dis\int_{\R^3} f(\eta(|\var_{n}z|)w(z))\eta(|\var_{n}z|)w(z)\rig\int_{\R^{3}}f(w)w
$$
and
\begin{equation}
\dis\int_{\R^3}
F(\eta(|\var_{n}z|)w(z))\rig\dis\int_{\R^{3}}F(w).\label{equation 8}
\end{equation}
We see that there is a subsequence of $\{t_{n}\}$ with $t_{\var_{n}}\rig 1$. In fact, since $\eta=1$ in $B_{\frac{\delta}{2}}(0)$ and
$B_{\frac{\delta}{2}}(0)\subset B_{\frac{\delta}{2\var_{n}}}(0)$ for $n\in \N$ large enough, it follows from $(\ref{equation 5})$ that
$$
\frac{M(t^{2}_{\var_{n}}\Lambda_{n}^{2})}{t^{2}_{\var_{n}}\Lambda_{n}^{2}}\geq \frac{1}{\Lambda_{n}^{4}}
\dis\int_{B_{\frac{\delta}{2}}(0)} \left[\frac{f(t_{\var_{n}}w(z))}{\left(t_{\var_{n}}w(z)\right)^{3}}\right]w(z)^{4},
$$
of continuity of $w$ (follows from standard regularity theory), there
is $\cha{z}\in \R^{3}$ such that
$w(\cha{z})=\dis\min_{\overline{B_{\frac{\delta}{2}}(0)}}w(z)$. So,
from $(f_{4})$
\begin{equation}
\frac{1}
{\Lambda_{n}^{4}}\frac{f(t_{\var_{n}}w(\cha{z}))}{\left(t_{\var_{n}}w(\cha{z})\right)^{3}}
\dis\int_{B_{\frac{\delta}{2}}(0)}w(z)^{4}\leq
\frac{M(t^{2}_{\var_{n}}\Lambda_{n}^{2})}{t^{2}_{\var_{n}}\Lambda_{n}^{2}}.\label{equation
6}
\end{equation}
Suppose by contradiction that there is a subsequence
$\{t_{\var_{n}}\}$ with $t_{\var_{n}}\rig \infty$. Thus, passing
to the limit as $n\rig \infty$ in $(\ref{equation 6})$, we conclude,
from $(M_{3})$ and $(f_{3})$, that  the left side converges to
infinity and the right side is bounded, which is a contradiction.
Therefore, $\{t_{\var_{n}}\}$ is bounded and passing a subsequence
we have $t_{\var_{n}}\rig t_{0}$ with $t_{0}\geq 0$.

From (\ref{equation 5}), (\ref{equation 4}), $(M_{1})$ and $(f_{4})$
we have that $t_{0}>0$.  Thus, passing to the limit as $n\rig
\infty$ in $(\ref{equation 5})$, we have
\begin{equation}
M(t_{0}^{2}\|w\|_{0}^{2})\|w\|_{0}^{2}t_{0}=\dis\int_{\R^{3}}f(t_{0}w)w.
\end{equation}
Since $w\in \mathcal{N}_{0}$, we obtain $t_{0}=1$. So, passing to
the limit of $n\rig \infty$ in $(\ref{equation 7})$ and using
$(\ref{equation 4})$ and $(\ref{equation 8})$ we obtain
$$
\lim_{n\rig
\infty}J_{\var_{n}}(\Upsilon_{\var_{n}}(y_{n}))=I_{0}(w)=c_{0},
$$
which is a contradiction with $(\ref{equation 2})$. $\fim$

Let us consider the following subset of the Nehari manifold

$$
\widetilde{\mathcal{N}}_{\var}=\{u\in \mathcal{N}_{\var}:
J_{\var}(u)\leq c_{0}+h_{1}(\var)\},
$$
where $h_{1}:\R_{+}\rig\R_{+}$ is a function such that $\Upsilon_{\var}(\Pi)\subset
\widetilde{\mathcal{N}}_{\var}$ and $\dis\lim_{\var\to 0}h_{1}(\var)=0$.
Observe that $h_{1}$ exists due to the Lemma \ref{Lemma 5.0.1}. In particular,
$\widetilde{\mathcal{N}}_{\var}\neq\emptyset$ for all small $\var>0$.

Now we consider $\rho > 0$ such that $\Pi_{\delta} \subset
B_{\rho}(0)$ and $\chi :  \R^3 \longrightarrow \R^3 $ defined by
\begin{eqnarray*}
\chi(x) =\ \  \left\{
             \begin{array}{l}
              x \quad se \quad |x| \leq \rho
        \\
        \\
        \dis\frac{\rho x}{|x|} \quad se \quad |x| \geq \rho.
        \\
        \\
             \end{array}
           \right.
\end{eqnarray*}

\noindent We also consider the barycenter map $\beta_{\var} :
\mathcal{N}_{\var} \longrightarrow \R^3 $ given by
\begin{eqnarray*}
\beta_{\var}(u) = \dis\frac{\dis\int_{\R^3}\chi(\var
x)u(x)^{2}}{\dis\int_{\R^3}u(x)^{2}}.
\end{eqnarray*}

Since $\Pi \subset B_{\rho}(0)$, the definition of $\chi$ and
Lebesgue's theorem imply that
\begin{eqnarray}\label{Lemma 5.0.2}
\dis\lim_{\var \rig 0}\beta_{\var}(\Upsilon_{\var}(y)) =y \
\mbox{uniformly in} \ y \in \Pi.
\end{eqnarray}

The next result is fundamental for showing that the solutions of
the auxiliary problem are solutions of the original problem.
Moreover, it allows us to show the behavior of such solutions.

\begin{Prop}\label{Proposition 5.0.1}
Let $\{u_{n}\}$ be a sequence in $H^{1}(\R^{3})$ such that
\begin{eqnarray*}
J_{\var_{n}}(u_{n})\rig c_{0}
\end{eqnarray*}
and
\begin{eqnarray*}
J'_{\var_{n}}(u_{n})(u_{n}) = 0, \ \forall n\in \N
\end{eqnarray*}
with $\var_{n} \rig 0$ when $n \rig \infty$. Then, there is a subsequence $\{\widetilde{y}_{n}\} \subset \R^{3}$
such that the sequence $v_{n}(x)= u_{n}(x + \tilde{y}_{n})$ has a convergent subsequence in $H^{1}(\R^{3})$. Moreover, passing to a subsequence,
\begin{eqnarray*}
y_{n} \rightarrow \widetilde{y} \ \mbox{with} \ y \in \Pi ,
\end{eqnarray*}
where $y_{n} = \var_{n}\widetilde{y}_{n}$.

\end{Prop}

\noindent {\bf Proof.} We can always to consider $u_{n}\geq 0$ and $u_{n}\neq 0$.
As in Lemma \ref{Lemma 4.2.1} and arguing as
Remark \ref{limitefraconaonulo} we have that $\{u_{n}\}$ is bounded
in $H^{1}(\R^{3})$ and there are $(\widetilde{y}_{n}) \subset \R^3$
and positive constants $R$ and $\alpha$ such that
\begin{eqnarray}\label{afirmation}
 \dis\liminf_{n\rig\infty}\dis\int_{B_{R}(\widetilde{y}_n)}u_{n}^{2} \geq \alpha >
 0.
\end{eqnarray}

Considering $v_{n}(x)=u_{n}(x+\widetilde{y}_{n})$ we conclude that
$\{v_{n}\}$ is bounded in $H^{1}(\R^{3})$ and therefore, passing to
a subsequence, we get
$$
v_{n}\righ v, \ \mbox{in} \ H^{1}(\R^{3})
$$
with $v\neq 0$. For each $n\in \N$, let $t_{n}>0$ such that
$\widetilde{v}_{n}=t_{n}v_{n}\in \mathcal{N}_{0}$ (see Lemma $\ref{Lemma 3.1.1} (A_{1})$). We have that
\begin{eqnarray*}
c_{0} &\leq& I_{0}(\widetilde{v}_{n})= \dis\frac{1}{2}\cha{M}(t^{2}_{n}\|u_{n}\|_{0}^{2})-\dis\int_{\R^3}F(t_{n}u_{n})\\
&\leq&
\dis\frac{1}{2}\cha{M}(t^{2}_{n}\|u_{n}\|_{\var_{n}}^{2})-\dis\int_{\R^3}G(\var_{n}
x, t_{n}u_{n}).
\end{eqnarray*}
Hence,
\begin{eqnarray}\label{vouusar}
c_{0} \leq I_{0}(\widetilde{v}_{n}) \leq J_{\var}(t_{n}u_{n}) \leq
J_{\var}(u_{n}) = c_{0} + o_{n}(1),
\end{eqnarray}
which implies,
\begin{equation}\label{convergence}
I_{0}(\widetilde{v}_{n}) \rig c_{0} \ \mbox{and} \
\{\widetilde{v}_{n}\} \subset \mathcal{N}_{0}.
\end{equation}

Thus, $\{\widetilde{v}_{n}\}$ is bounded in $H^{1}(\R^{3})$ and
$\widetilde{v}_{n}\righ \widetilde{v}$. From well-known arguments we
can assume that $t_{n}\rig t_{0}$ with $t_{0}>0$. So, from
uniqueness of the weak limit we have $\widetilde{v}=t_{0}v$, $v\neq
0$. From Lemma $\ref{compacity Lemma}$ we obtain,
\begin{equation}
\widetilde{v}_{n}\rig \widetilde{v} \ \mbox{in} \
H^{1}(\R^{3}).\label{convergence 1}
\end{equation}

This convergence implies
$$
v_{n}\rig \frac{\widetilde{v}}{t_{0}}=v \ \mbox{in} \ H^{1}(\R^{3})
$$
and
\begin{equation}
I_{0}(\widetilde{v})=c_{0} \ \mbox{and} \
I'_{0}(\widetilde{v})\widetilde{v}=0.
\end{equation}

Now, we will show that $\{y_{n}\}$ is bounded, where
$y_{n}=\var_{n}\widetilde{y}_{n}$. In fact, otherwise, there exists
a subsequence $\{y_{n}\}$ with $|y_{n}|\rig \infty$. Observe that
$$
m_{0}\|v_{n}\|_{0}^{2}\leq \int_{\R^{3}}g(\var_{n}z+y_{n},
v_{n})v_{n}.
$$
Let $R>0$ such that $\Omega\subset B_{R}(0)$. Since we may suppose
that $|y_{n}|\geq 2R$, for each $z\in B_{\frac{R}{\var_{n}}}(0)$ we
have
$$
|\var_{n}z+y_{n}|\geq |y_{n}|-|\var_{n}z|\geq 2R-R=R.
$$
Thus,
$$
m_{0}\|v_{n}\|_{0}^{2}\leq
\int_{B_{\frac{R}{\var_{n}}}(0)}\widetilde{f}(v_{n})v_{n}+\dis\int_{\R^{3}\backslash
B_{\frac{R}{\var_{n}}}(0)}f(v_{n})v_{n}.
$$
Since $v_{n}\rig v$ in $H^{1}(\R^{3})$, it follows from  Lebesgue's
theorem that
$$
\dis\int_{\R^{3}\backslash
B_{\frac{R}{\var_{n}}}(0)}f(v_{n})v_{n}=o_{n}(1).
$$
On the other hand, since $\widetilde{f}(v_{n})\leq
\frac{V_{0}}{K}v_{n}$, we obtain
$$
m_{0}\|v_{n}\|_{0}^{2}\leq
\frac{1}{K}\int_{B_{\frac{R}{\var_{n}}}(0)}V_{0}v_{n}^{2}+o_{n}(1),
$$
and therefore,
$$
\left(m_{0}-\frac{1}{K}\right)\|v_{n}\|_{0}\leq o_{n}(1),
$$
which is a contradiction. Hence, $\{y_{n}\}$ is bounded and we can
assume $y_{n}\rig \overline{y}$ in $\R^{3}$. We see that
$\overline{y}\in\overline{\Omega}$ because if
$\overline{y}\notin\overline{\Omega}$, we can proceed as above and
conclude that  $\|v_{n}\|_{0}\leq o_{n}(1)$.

In order to prove that $V(\overline{y})=V_{0}$, we suppose by
contradiction that $V_{0}<V(\overline{y})$. Consequently, from
$(\ref{convergence 1})$, Fatou's Lemma and the invariance of
$\mathbb{R}^{3}$ by translations, we obtain
\begin{eqnarray*}
c_{0}&<&\dis\liminf_{n\to \infty} \left[\frac{1}{2}\cha{M}\left(
\int_{\R^{3}}|\nabla
\widetilde{v}_{n}|^{2}+\int_{\R^{3}}V(\var_{n}z+y_{n})
\widetilde{v}_{n}^{2}\right)
-\dis\int_{\R^{3}}F(\widetilde{v}_{n})\right]\\
&\leq&\liminf_{n\rig \infty}J_{\var_{n}}(t_{n}u_{n})\\
&\leq& \liminf_{n\rig \infty}J_{\var_{n}}(u_{n})=c_{0},
\end{eqnarray*}
which is a contradiction and the proof is finished. $\fim$

\begin{Co}\label{KKK} Assume the same hypotheses of
Proposition \ref{Proposition 5.0.1}. Then, for any given $\gamma_{2}
>0$, there exists $R>0$ and $n_{0}\in \mathbb{N}$ such that
$$
\dis\int_{B_R(\widetilde{y}_n)^c}\left( |\nabla u_n|^2 + |u_n|^2
\right) <\gamma_{2},~~\mbox{for all } n\geq n_{0}.
$$
\end{Co}
\noindent {\bf Proof.} By using the same notation of the proof of
Proposition \ref{Proposition 5.0.1}, we have for any $R>0$
\begin{eqnarray*}
\dis\int_{B_R(\widetilde{y}_n)^c}\left( |\nabla u_n|^2 + |u_n|^2
\right) & = & \dis\int_{B_R(0)^c} (|\nabla v_n|^2 + |v_n|^2).
\end{eqnarray*}
Since $(v_n)$ strongly converges in $H^1(\mathbb{R}^N)$ the result
follows.$\fim$

\begin{Le}\label{Lemma 5.0.3}
Let $\delta > 0$ and $\Pi_{\delta}=\{x \in \R^3 : dist(x,M) \leq
\delta\}$. Then,
\begin{eqnarray*}
\dis\lim_{\var \rightarrow 0}\dis\sup_{u \in
\mathcal{\widetilde{N}}_{\var}}\dis\inf_{y \in
\Pi_{\delta}}|\beta_{\var}(u) - y| = \dis\lim_{\var \rightarrow 0}\dis\sup_{u \in
\mathcal{\widetilde{N}}_{\var}}dist(\beta_{\var}(u),\Pi_{\delta}) = 0.
\end{eqnarray*}
\end{Le}

\noindent {\bf Proof.} The proof of this Lemma follows from well-known
arguments and can be found in \cite[Lemma 3.7]{Alves and Figueiredo
1}. $\fim$

\subsection{Multiplicity of solutions for $(P_{\var,A})$}

In the sequel we prove our multiplicity result for the problem
$(P_{\var,A})$, by using arguments slightly different to those in
$\cite{Wang 1}$, in fact, since $S_{\var}^{+}$ is a incomplete metric space,
we can't use (directly) an abstract result as in [$\cite{Cingolani}$, Theorem 2.1], instead,
we invoke the category abstract result in $\cite{Szulkin}$.

\begin{Th}\label{Theorem 5.0.1}
Suppose that the function $M$ satisfies $(M_{1})-(M_{3})$, the
potential $V$ satisfies $(V_{1})-(V_{2})$ and the function $f$
satisfies $(f_{1})-(f_{4})$. Then, given $\delta>0$ there is
$\overline{\var}=\overline{\var}(\delta)>0$ such that the auxiliary
problem $(P_{\var,A})$ has at least $Cat_{\Pi_{\delta}}(\Pi)$
positive solutions, for all $\var\in (0,\overline{\var})$.
\end{Th}

\noindent {\bf Proof.}  For each $\var>0$, we define the function $\zeta_{\var}:\Pi\rightarrow S_{\var}^{+}$ by
$$
\zeta_{\var}(y)=m_{\var}^{-1}(\Upsilon_{\var}(y)),
\ \forall y\in \Pi.
$$
From the Lemma $\ref{Lemma 5.0.1}$, we have
$$
\lim_{\var\to 0}\Psi_{\var}(\zeta_{\var}(y))=\lim_{\var\to
0}J_{\var}(\Upsilon_{\var}(y))=c_{0}, \ \mbox{uniformly in} \ y\in\Pi.
$$
Thus, the set
$$
\widetilde{S}_{\var}^{+}=\{u\in S_{\var}^{+}: \Psi_{\var}(u)\leq c_{0}+
h_{1}(\var)\},
$$
is nonempty, for all $\var\in (0, \overline{\var})$, because $\zeta_{\var}(\Pi)\subset \widetilde{S}_{\var}^{+}$,
where the function $h_{1}$ was already introduced in the definition of the set $\widetilde{\mathcal{N}}_{\var}$.

From above considerations, together with
Lemma $\ref{Lemma 5.0.1}$,  Lemma $\ref{lema3.2.2}(A_{3})$, equality $(\ref{Lemma 5.0.2})$ and Lemma $\ref{Lemma 5.0.3}$, there is
$\overline{\var}=\overline{\var}(\delta)>0$, such that the diagram of continuous applications bellow is well defined for $\var\in(0, \overline{\var})$
$$
\Pi
\stackrel{\Upsilon_{\var}}{\longrightarrow}\Upsilon_{\var}(\Pi)\stackrel{{m_{\var}^{-1}}}{\longrightarrow}\zeta_{\var}(\Pi)\stackrel{{m_{\var}}}{\longrightarrow}
\Upsilon_{\var}(\Pi)\stackrel{{\beta_{\var}}}{\longrightarrow}\Pi_{\delta}.
$$

We conclude from ($\ref{Lemma 5.0.2}$) that there is a function $\lambda(\var, y)$ with
$|\lambda(\var, y)|<\frac{\delta}{2}$ uniformly in $y\in \Pi$, for all $\var\in(0, \overline{\var})$,
such that $\beta_{\var}(\Upsilon_{\var}(y))=y+\lambda(\var, y)$ for all $y\in \Pi$. Hence,
the application $H:[0,1]\times \Pi\rig \Pi_{\delta}$ defined by
$H(t,y)=y+(1-t)\lambda(\var, y)$ is a homotopy between
$\alpha_{\var}\circ\zeta_{\var}=\beta_{\var}\circ\Upsilon_{\var}$
and the inclusion $i: \Pi\rig \Pi_{\delta}$, where $\alpha_{\var}=\beta_{\var}\circ m_{\var}$. Therefore,
\begin{equation}
cat_{\zeta_{\var}(\Pi)}\zeta_{\var}(\Pi)\geq cat_{\Pi_{\delta}}(\Pi).\label{ineq 4}
\end{equation}

It follows from Corollary $\ref{Corollary 4.2.1}$ and from category
abstract theorem in \cite{Szulkin}, with $c=c_{\var}\leq c_{0}+h_{1}(\var)=d$
and $K=\zeta_{\var}(\Pi)$, that $\Psi_{\var}$ has at least
$cat_{\zeta_{\var}(\Pi)}\zeta_{\var}(\Pi)$ critical points on
$\widetilde{S}_{\var}^{+}$. So, from item $(d)$ of the Proposition
$\ref{Proposition 4.1.2}$ and from $(\ref{ineq 4})$, we conclude
that $J_{\var}$ has at least $cat_{\Pi_{\delta}}(\Pi)$ critical
points in $\widetilde{\mathcal{N}}_{\var}$.

%On the other hand, suppose that $cat_{\Pi_{\delta}}(\Pi)>1$. For each $\var\in (0, \overline{\var})$,
%follows from Lemma $\ref{lema3.2.2}(A_{4})$ that there are $u_{\var}\in S_{\var}^{+}$  and
%$h_{2}:\R_{+}\rig \R_{+}$ such that $\lim_{\var\to 0}h_{2}(\var)=0$ and
%$$
%c_{0}+h_{1}(\var)< \Psi_{\var}(u_{\var})\leq c_{0}+h_{2}(\var).
%$$
%Denoting,
%$$
%h_{2}(\var)\},
%$$
%we have $\widetilde{S}_{\var}^{+}\varsubsetneq\cha{S}_{\var}^{+}$. Moreover, we define the
%homotopy $\eta_{\var}:[0,1]\times K\rig S_{\var}^{+}$ by
%$$
%\eta_{\var}(t,u)=\frac{t u_{\var}+(1-t) u}{\| t u_{\var}+(1-t) u\|_{\var}}.
%$$
%Note that $\eta_{\var}(0, u)=u$ and $\eta_{\var}(1, u)=u_{\var}$, for all $u\in K$. Moreover,

$\fim$

\section{Proof of Theorem \ref{Theorem 1.1}}

In this section we prove our main theorem. The idea is to show that
the solutions obtained in Theorem \ref{Theorem 5.0.1} verify the
following estimate $u_{\var}(x) \leq a \,\, \forall x \in
\Omega^{c}_{\var}$ for $\var$ small enough. This fact implies that
these solutions are in fact solutions of the original problem
$(\widetilde{P}_{\var})$. The key ingredient is the following
result, whose proof uses an adaptation of the arguments found in
\cite{Gongbao}, which are related to the Moser iteration method
\cite{Moser} .

\begin{Le}\label{leminha1} Let $\var_n \to 0^+$ and $u_{n} \in \widetilde{\mathcal{N}}_{\var_n}$ be a solution of $(P_{\var_{n},A})$.
Then $J_{\var_n}(u_n) \to c_0$ and $u_n \in
L^{\infty}(\mathbb{R}^3)$. Moreover, for any given $\gamma >0$,
there exists $R>0$ and $n_{0}\in \mathbb{N}$ such that
\begin{equation} \label{eq1_leminha1}
|u_{n}|_{L^{\infty}(B_{R}(\tilde{y}_{n})^c)}<\gamma,~~~~\mbox{for
all} \ n\geq n_{0},
\end{equation}
where $\tilde{y}_{n}$ is given by Proposition \ref{Proposition
5.0.1}.
\end{Le}
\noindent {\bf Proof.} Since $J_{\var_n}(u_n) \leq c_0 + h(\var_n)$
with $\dis\lim_{n\to\infty} h(\var_n)= 0$, we can argue as in the
proof of the inequality (\ref{vouusar}) to conclude that $J_{\var_n}(u_n)
\to c_0$. Thus, we may invoke  Proposition \ref{Proposition 5.0.1}
to obtain a sequence $(\widetilde{y}_n) \subset \mathbb{R}^3$
satisfying the conclusions of that proposition.

Fix $R>1$ and consider $\eta_{R} \in C^{\infty}(\mathbb{R}^3)$ such
that $0 \leq \eta_R \leq 1$, $\eta_{R} \equiv  0$ in $B_{R/2}(0)$,
$\eta_R \equiv 1$ in $B_{R}(0)^c$ and  $|\nabla \eta_{R}|\leq C/R$.
For each $n \in\mathbb{N}$ and $L>0$, we define
$\eta_{n}(x):=\eta_{R}(x-\tilde{y}_{n})$, $u_{L,n} \in
H^1(\mathbb{R}^3)$ and $z_{L,n} \in H_{\var}$ by setting
$$
u_{L,n}(x) := \min\{ u_n(x), L \},~~~~
 z_{L,n} :=
\eta_{n}^{2}u_{L,n}^{2(\beta - 1)}u_{n},
$$
with $\beta > 1$ to be determined later.

From definition of $z_{L,n}$ and $J'_{\var_{n}}(u_{n})z_{L,n}=0$,
we have
\begin{equation*}
\begin{array}{c}
m_{0}\biggl[\dis\int_{\mathbb{R}^{3}} \eta_{n}^{2}u_{L,n}^{2(\beta -
1)}|\nabla u_{n}|^{2} +
 2 \dis\int_{\mathbb{R}^{3}} \eta_{n}u_{n}u_{L,n}^{2(\beta - 1)} \nabla
\eta_{n}\cdot\nabla u_{n}\biggl] \ \  \\
\leq \dis\int_{\mathbb{R}^{3}} \left( g(\var_n
x,u_{n})-m_{0}V(\var_{n}x)u_{n}
\right)\eta_{n}^{2}u_{n}u_{L,n}^{2(\beta - 1)} .
\end{array}
\end{equation*}

Now, the result follows arguing as in \cite[Lemma 4.1]{CPDE}.$\fim$

We are now ready to prove the main result of the paper.

\subsection{Proof of Theorem \ref{Theorem 1.1}}

\noindent Suppose that $\delta>0$ is such that $\Pi_{\delta} \subset
\Omega$. We first claim that there exists
$\widetilde{\var}_{\delta}>0$ such that, for any
$0<\var<\widetilde{\var}_{\delta}$ and any solution $u \in
\widetilde{\mathcal{N}}_{\var}$ of the problem $(P_{\var,A})$,
there holds
\begin{equation} \label{prova1}
| u |_{L^{\infty}(\mathbb{R}^3 \setminus \Omega_{\var})} < a.
\end{equation}
In order to prove the claim we argue by contradiction. So, suppose
that for some sequence $\var_n \to 0^+$ we can obtain $u_n \in
\widetilde{\mathcal{N}}_{\var_n}$ such that $J_{\var_n}'(u_n) =0$
and
\begin{equation}
\label{provafinal_eq1} | u_n |_{L^{\infty}(\mathbb{R}^3 \setminus
\Omega_{\var_n})} \geq a.
\end{equation}
As in Lemma \ref{leminha1}, we have that $J_{\var_n}(u_n) \to c_0$
and therefore we can use Proposition \ref{Proposition 5.0.1} to
obtain a sequence $(\widetilde{y}_n) \subset \mathbb{R}^3$ such that
$\var_n \widetilde{y}_n \to y_0 \in \Pi$.

If we take $r>0$ such that $B_r(y_0) \subset B_{2r}(y_0) \subset
\Omega$ we have that
$$
B_{r/\var_n}(y_{0}/\var_n)=\frac{1}{\var_n}B_{r}(y_0)\subset
\Omega_{\var_n}.
$$
Moreover, for any $z \in B_{r/\var_n}(\widetilde{y}_{n})$, there
holds
$$
\left|z-\frac{y_0}{\var_{n}}\right|\leq
|z-\widetilde{y}_{n}|+\left|\tilde{y}_{n}-\frac{y_0}{\var_{n}}\right|<
\frac{1}{\var_{n}}(r+o_{n}(1))< \frac{2r}{\var_{n}},
$$
for $n$ large. For these values of $n$ we have that
$B_{r/\var_n}(\widetilde{y}_n) \subset \Omega_{\var_n}$ or,
equivalently, $\mathbb{R}^3 \setminus \Omega_{\var_n} \subset
\mathbb{R}^3 \setminus B_{r/\var_n}(\widetilde{y}_n)$. On the other
hand, it follows from Lemma \ref{leminha1} with $\gamma=a$ that, for
any $n \geq n_0$ such that $r/\var_n > R$, there holds
$$
|u_n |_{L^{\infty}(\mathbb{R}^3 \setminus \Omega_{\var_n})} \leq
|u_n |_{L^{\infty}({\mathbb{R}^3 \setminus
B_{r/\var_n}(\widetilde{y}_n)})} \leq |u_n
|_{L^{\infty}(\mathbb{R}^3 \setminus B_R(\widetilde{y}_n))} < a,
$$
which contradicts (\ref{provafinal_eq1}) and proves the claim.

Let $\widehat{\var}_{\delta}>0$ given by Theorem \ref{Theorem 5.0.1}
and set $\var_{\delta} := \min\{ \widehat{\var}_{\delta},
\widetilde{\var}_{\delta}\}$. We shall prove the theorem for this
choice of $\var_{\delta}$. Let $0<\var<\var_{\delta}$ be fixed. By
applying Theorem \ref{Theorem 5.0.1} we obtain
cat$_{\Pi_{\delta}}(\Pi)$ nontrivial solutions of the problem
$(P_{\var,A})$. If $u \in H_{\var}$ is one of these solutions we
have that $u \in \widetilde{\mathcal{N}}_{\var}$, and therefore we
can use (\ref{prova1}) and the definition of $g$ to conclude that
$g_{\var}(\cdot,u) \equiv f(u)$. Hence, $u$ is also a solution of
the problem $(\widetilde{P}_{\var})$. An easy calculation shows that
$\widehat{u}(x):=u(x/\var)$ is a solution of the original problem
$(P_{\var})$. Then, $(P_{\var})$ has at least
cat$_{\Pi_{\delta}}(\Pi)$ nontrivial solutions.

We now consider $\var_n \to 0^+$ and take a sequence $u_n \in
H_{\var_n}$ of solutions of the problem $(\widetilde{P}_{\var_n})$
as above. In order to study the behavior of the maximum points of
$u_{n}$, we first notice that, by $(g_1)$, there exists $\gamma>0$
such that
\begin{equation} \label{comp_maximo1}
g(\var x,s)s \leq \frac{V_0}{K}s^{2},~~\mbox{for all } x
\in\mathbb{R}^3,\,s \leq \gamma.
\end{equation}
By applying Lemma \ref{leminha1} we obtain $R>0$ and
$(\widetilde{y}_n) \subset \mathbb{R}^3$ such that
\begin{equation} \label{comp_maximo2}
| u_n|_{L^{\infty}(B_R(\widetilde{y}_n))^c} < \gamma,
\end{equation}
Up to a subsequence, we may also assume that
\begin{equation} \label{comp_maximo3}
|u_n |_{L^{\infty}(B_R(\widetilde{y}_n))} \geq \gamma.
\end{equation}
Indeed, if this is not the case, we have
$|u_n|_{L^{\infty}(\mathbb{R}^3)} < \gamma$, and therefore it
follows from $J_{\var_n}'(u_n)=0$ and (\ref{comp_maximo1})  that
$$
m_{0} \| u_n\|_{\var_n}^2 \leq \int_{\displaystyle_{\mathbb{R}^{3}}}
g(\var_{n}x,u_n)u_n \leq \frac{V_0}{K}
\int_{\displaystyle_{\mathbb{R}^{3}}} u_n^2.
$$
The above expression implies that $\|u_n\|_{\var_{n}}=0$, which does
not make sense. Thus, (\ref{comp_maximo3}) holds.

By using (\ref{comp_maximo2}) and (\ref{comp_maximo3}) we conclude
that the maximum point $p_n \in \mathbb{R}^3$ of $u_{n}$ belongs to
$B_R(\widetilde{y}_n)$. Hence $p_n = \widetilde{y}_n + q_n$, for
some $q_n \in B_R(0)$. Recalling that the associated solution of
$(P_{\var_n})$ is of the form $\widehat{u}_n(x)=u_n(x/\var_n)$, we
conclude that the maximum point $\eta_{n}$ of $\widehat{u}_n$ is
$\eta_n := \var_n \widetilde{y}_n + \var_n q_n$. Since $(q_n)
\subset B_R(0)$ is bounded and $\var_n \widetilde{y}_n \to y_0 \in
\Pi$ (according to Proposition \ref{Proposition 5.0.1}), we obtain
$$
\dis\lim_{n\rightarrow \infty}V(\eta_{\var_{n}})=V(y_0)=V_{0},
$$
which concludes the proof of the theorem.  $\fim$

\vspace{0.5cm} \noindent {\bf Acknowledgement.} This work was done
while the authors were visiting the "Departamento de ecuaciones
diferenciales y an\'{a}lisis num\'{e}rico" of the "Universidad de
Sevilla". They would like to express their gratitude to the Prof.
Antonio Suarez for his warm hospitality.

\end{document}